\documentclass[a4paper]{article}
\usepackage[utf8]{inputenc}
\usepackage{amsmath,amsthm,amssymb,amsfonts}

\newtheorem{theorem}{Theorem}
\newtheorem{lemma}{Lemma}
\newtheorem{claim}{Claim}
\newtheorem{corol}{Corollary}

\usepackage{capt-of}

\usepackage{tikz} 
\usetikzlibrary{positioning}

\newcommand{\Exp}{{\sf E}}
\newcommand{\Var}{\mathrm{Var}}

\usepackage{color}

\title{Weak saturation stability}
\author{O. Kalinichenko\footnote{Moscow Institute of Physics and Technology, Dolgoprudny, Russia}, M. Zhukovskii\footnotemark[\value{footnote}] \footnote{Adyghe State University, Caucasus mathematical center, Maykop, Republic of Adygea, Russia} \footnote{The Russian Presidential Academy of National Economy and Public Administration, Moscow, Russia}}
\date{}

\begin{document}

\maketitle

\begin{abstract}
 %  We study wsat$(G,H)$, the minimum number of edges in a weakly $H$-saturated subgraph of $G$. We prove that wsat$(K_n,H)$ is stable (remains the same after independent removal of every edge of $K_n$ with a constant probability) for all pattern graphs $H$ such that there exists a `local' set of edges that percolates in $K_n$ (this is true, for example, for cliques and complete bipartite graphs). Also, we find a threshold probability for the weak $K_{1,t}$-saturation stability. 
    
    The paper studies wsat$(G,H)$ which is the minimum number of edges in a weakly $H$-saturated subgraph of $G$. We prove that wsat$(K_n,H)$ is `stable' --- remains the same after independent removal of every edge of $K_n$ with constant probability --- for all pattern graphs $H$  such that there exists a `local' set of edges percolating in $K_n$. This is true, for example, for cliques and complete bipartite graphs. We also find a threshold probability for the weak $K_{1,t}$-saturation stability. 

\end{abstract}

\vspace{0.5cm}

%Let $G$ and $H$ be graphs. A spanning subgraph $H \subset G$ is \textit{weakly $(G, H)$-saturated}, if it contains no copy of $H$, but there exists an ordering $e_1, ...e_m$ of the edges of $G \setminus F$ such that $\forall i \in {1, ...m} \ F \cup \{e_1, ...e_i\}$ contains a copy of $H$ that contains $e_i$. ... 

\section{Introduction}

Let $G$ and $H$ be graphs (below, we refer to them as {\it host} and {\it pattern} graphs respectively). Let $F\subset G$ be a spanning subgraph of $G$. Let us call a sequence of graphs $F=F_0\subset\ldots\subset F_m=G$ an {\it $H$-bootstrap percolation process}, if $F_i$ is obtained from $F_{i-1}$ by adding an edge that belongs to a copy of $H$ in $F_i$. $F$ is {\it weakly $(G, H)$-saturated}, if %it contains no copy of $H$ but
$G$ can be obtained from $F$ in an $H$-bootstrap percolation process  (i.e., there exists an ordering $e_1,\ldots,e_m$ of the edges of $G \setminus F$ such that, for every $i\in[m]:=\{1,\ldots,m\}$, $F \sqcup \{e_1,\ldots,e_i\}$ has a copy of $H$ that contains $e_i$).  The smallest number of edges in a weakly $(G, H)$-saturated graph is called {\it the weak saturation number} and is denoted by $\mathrm{wsat}(G, H)$. %We will also call $H$ a \textit{pattern}.
This notion was first introduced by Bollob\'{a}s in 1968 \cite{bollobas}.\\

In this paper, we study the phenomena of {\it stability} of the weak saturation number. It was observed by Kor\'{a}ndi and Sudakov~\cite{sudakov} that $\mathrm{wsat}(G=K_n,K_s)$ remains the same after the deletion of each edge of $K_n$ independently with a constant positive probability (as usual, $K_n$ denotes a complete graph on $n$ vertices). We show that this stability property holds for a wider class of pattern graphs $H$, and conjecture that it actually holds for all $H$. We also find a threshold probability for the stability of the weak $(K_n,K_{1,t})$-saturation number. Before stating the results, let us recall known values of the weak saturation number when $G=K_n$. We denote by $K_{s,t}$ a complete bipartite graph with parts of size $s$ and $t$.\\

The exact value of $\mathrm{wsat}(K_n, K_s)$ was achieved by Lov\'{a}sz \cite{lovasz}% using algebraic techniques
: if $n \ge s \ge 2$ then 
$$
\mathrm{wsat}(K_n, K_s) = {n\choose 2} - {n-s+2\choose 2}.
$$

%For bipartite graphs, a related problem of \textit{bisaturation} has also been well studied. Let $H = (V_1\sqcup V_2, E)$ be a bipartite graph. A spanning subgraph $F$ of $K_{n,\ell}$ (complete bipartite graph with parts of size $n$ and $\ell$) is \textit{weakly $(K_{n, \ell}, H)$-bisaturated} if there exists an ordering $e_1,\ldots,e_m$ of the edges of $K_{n, \ell} \setminus F$ such that, for every $i\in[m]$, $F \sqcup \{e_1, \ldots,e_i\}$ contains a copy $H' = (V_1' \sqcup V_2', E')$ of $H$ that contains $e_i$ and such that $V_1'$ and $V_2'$ lie in the first and second classes of the vertices of $K_{n, \ell}$ respectively. 

%The \textit{bisaturation number}, denoted by $w(K_{n, \ell}, H)$, is the minimum number of edges in a weakly $(K_{n, \ell}, H)$-bisaturated graph. It is obvious that a weakly $(K_{n, \ell}, H)$-bisaturated graph is weakly $(K_{n, \ell}, H)$-saturated which implies that $w(K_{n, \ell}, H) \ge \mathrm{wsat}(K_{n, \ell}, H)$. Alon proved \cite{alon} that 
%$$
%w(K_{n, \ell}, K_{s, t}) = n\cdot \ell - (n - s + 1)(\ell - t + 1)
%$$ 
%for $2 \le s \le t$ and $2 \le n \le \ell$.  Moshkovitz and Shapira \cite{moshkovitz} deduced from this result that 
%$$
%\mathrm{wsat}(K_{n,n}, K_{s,t}) = n^2 - (n - s + 1)^2 + (t - s)^2
%$$ 
%for $2 \le s \le t \le n$. %?? ????? ??? ????????????
%This result can be easily generalized (see Appendix A in \cite{martins})

%$$
%\mathrm{wsat}(K_{n,\ell}, K_{s,t}) = (n +\ell - s + 1)(s-1)  + (t - s)^2
%$$ 
%for $2 \le s \le t , n \le \ell$.

The value of $\mathrm{wsat}(K_n,K_{s,t})$ for an arbitrary choice of parameters is still unknown. The most general result was obtained by Kalai~\cite{Kalai} in 1985 and Kronenberg, Martins and Morrison \cite{martins} in 2020. They proved that 
$$
\mathrm{wsat}(K_n, K_{t,t}) = (t - 1)(n + 1 - t/2),
$$ 
$$
\mathrm{wsat}(K_n, K_{t,t+1}) = (t - 1)(n + 1 - t/2) + 1
$$ 
%if $t \ge 2$ and $n \ge 3t - 3$. They also proved a similar result for $K_{t,t+1}$: 
%$$
%\mathrm{wsat}(K_n, K_{t,t+1}) = (t - 1)(n + 1 - t/2) + 1
%$$ 
if $t \ge 2$ and $n \ge 3t - 3$. In~\cite{martins} general bounds for arbitrary choice of parameters $s,t$ were also obtained:
\begin{equation}
\mathrm{wsat}(K_n, K_{s,t}) \le (s - 1)(n - s) + {t \choose 2}
\label{general_bipartite_upper}
\end{equation}
if $t > s \ge 2$ and $n \ge 2(s + t) - 3$ and
\begin{equation}
\mathrm{wsat}(K_n, K_{s,t}) \ge (s - 1)(n - t + 1) + {t \choose 2}
\label{general_bipartite_lower}
\end{equation} 
if  $t > s \ge 2$ and $n \ge 3t - 3$.

Notice that, for $s=1$, the exact value of the weak saturation number is straightforward: 
$$
\mathrm{wsat}(K_n, K_{1, t}) = {t\choose2}.
$$ 
Indeed, $\mathrm{wsat} (K_n, K_{1,t}) \le {t \choose 2}$ since
we can use $F=K_t$ to restore all edges of $K_n$ (first, we restore all edges with one endpoint in the clique, then all the rest). Eventually, $\mathrm{wsat} (K_n, K_{1,t}) \ge {t \choose 2}$ since any weakly $(K_n, K_{1, t})$-saturated graph $F$ must contain vertices $v_1,\ldots,v_{t-1}$ such that, for every $i\in[t-1]$, the number of neighbors of $v_i$ in $V(F)\setminus\{v_1,\ldots,v_{i-1}\}$ is at least $t-i$. To see this, it is sufficient to consider the first distinct $t-1$ vertices playing the roles of the central vertices of $K_{1,t}$ in a bootstrap percolation process that starts on $F$ and finishes on $K_n$.\\

%Let's call the only vertex in the left class of $K_{1, s}$ a \textit{central} vertex.

%It is obvious that the vertices that are central in the copies of $H$ that are formed during the adding of new edges must have degree at least $s-1$ at the time of adding. If $F$ is weakly $(K_n, K_{1, s})$-saturated, then the first $s-1$ unique central  vertices in the forming copies of $H$ can only be the vertices of $F$, because other vertices still have degree less than $s-1$. Also, they must have degree at least $s-1$ in $F$, because no edges can be added to a vertex of $F$ that has not yet been central (the vertices of the outside of $F$ can not yet be used as central ones). So, there must be at least $(s-1) + (s-2) + ... 1 = {s \choose 2}$ edges in $F$. It implies that $$\mathrm{wsat} (K_n, K_{1,s}) \ge {s \choose 2}$$

There are also many results about weak saturation numbers for other specific pairs of host and pattern graphs (e.g., for both $G$ and $H$ being complete bipartite~\cite{moshkovitz}, for multipartite graphs \cite{alon}, for disjoint copies of graphs \cite{faudree}, for hypercubes and grids \cite{balogh, pete, hypercube}). The weak saturation number for hypergraphs has also been studied \cite{alon, balogh, uniform, faudree, moshkovitz, tuza, tuza1}.\\

%For example, we know the values of the weak saturation numbers for complete multipartite graphs and hypergraphs [2], asymptotics of the weak saturation number of hypergraphs [29], the weak saturation number for families of hypergraphs with a fixed number of edges [12, 25, 30], complete bipartite hypergraphs (in complete bipartite hypergraphs) [6, 24], pyramids in hypergraphs [26], families of graphs in the complete graph [27, 10, 28], families of disjoint copies of graphs [15], and for the case where $H$ is a hypercube or a grid [6, 7, 23]. 
% ?????? ?? ??? ?????? ??? ??????????, ????????!!! ? ??? ??? ???????? ?????...
%$\,$\\

In 2017,  Kor\'{a}ndi  and Sudakov \cite{sudakov} proved that, if $s \ge 3$, then $\mathrm{wsat}(K_n,K_s)$ is {\it stable}, i.e., for constant $p\in(0,1)$,
$$
\mathrm{wsat}(G(n,p), K_s) = \mathrm{wsat}(K_n,K_s)
$$
with high probability.
Here, we denote by $G(n,p)$ the binomial random graph on the vertex set $[n]$, where every pair of distinct $i,j\in[n]$ is adjacent with probability $p$ independently of the others. Hereinafter, we say that some property holds {\it with high probability}, or {\it whp}, if its probability tends to 1 as $n\to\infty$.\\ %\textcolor{blue}{We call subgraph $F$ of $G(n, p)$ (which is a random element on the set of graphs with $n$ vertices) weakly $(G(n, p). H)$-saturated if for every realisation $G'$ of $G(n, p)$ the realisation $F'$ of $F$ is weakly $(G', H)$-saturated.}

%An attempt to research the value of $wsat(G, H)$ in the case of a random graph $G = G(n, p)$ was taken by Kor\'{a}ndi  and Sudakov \cite{6} in 2017. They found that for a constant $p \in (0, 1)$ w.h.p. $wsat(G(n,p), K_s) = wsat(K_n,K_s)$. 

In this paper, we prove a transference theorem %result
that can be used to derive such stability results. It immediately implies the result of Kor\'{a}ndi  and Sudakov as well as stability results for all complete bipartite pattern graphs (despite the fact that the exact value of $\mathrm{wsat}(K_n,K_{s,t})$ is not known for almost all pairs $s$ and $t$). Below, we denote by $\delta(H)$ the minimum degree of graph $H$. Without loss of generality, we set $V(K_n) = [n]$.\\

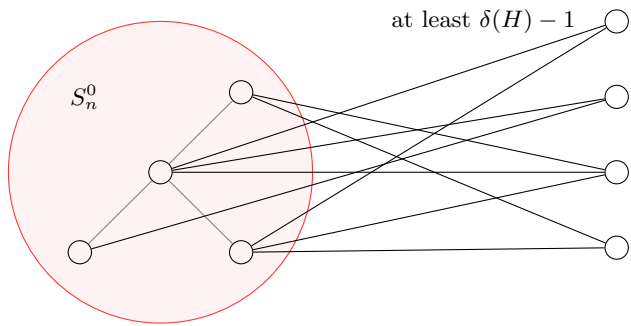
\begin{figure}
\begin{centering}
\begin{tikzpicture}[main/.style = {node distance={15mm}, draw, circle}] \small

\filldraw[color=red!80, fill=red!05] (0, 0) circle (2.0);

\node[main] (1) {}; 
\node[main] (2) [above right of=1] {}; 
\node[main] (3) [below left of=1]{}; 
\node[main] (4) [below right of=1]{}; 
\draw[gray] (1) -- (2);
\draw[gray] (1) -- (3);
\draw[gray] (1) -- (4);

\node[main] at (6, 0)   (a) {};
\node[main] at (6, 1)   (b) {};
\node[main] at (6, -1)  (c)     {};
\node[main] at (6, 2)  (d)     {};

\node[] at (-1, 1) (222) {$S^0_n$}; 

\draw[black] (a) -- (1);
\draw[black] (a)-- (2);
\draw[black] (a) -- (4);

\draw[black] (b) -- (1);
\draw[black] (b) {} -- (3);

\draw[black] (c) -- (2);
\draw[black] (c) -- (4);

\draw[black] (d) -- (1);
\draw[black] (d) -- (4);

\node[] at (4.25, 2) (222) {at least $\delta(H)-1$};

\end{tikzpicture} 
%\captionof{figure}{Structure of $F^0_n$.}
\label{fg:1}
\caption{structure of $F^0_n$.}
\end{centering}
\end{figure}

%\medskip

\begin{theorem}%[Transference]
Let $H$ be a graph without isolated vertices, and let $p \in (0,1)$, $C\geq\delta(H)-1$ be constants. For every $n\in\mathbb{N}$, let $F^0_n$ be a weakly $(K_n, H)$-saturated graph containing a set of vertices $S^0_n\subset[n]$ with $|S^0_n| \leq C$, such that every vertex from $[n] \setminus S^0_n$ is adjacent to  at least $\delta(H) - 1$ vertices of $S^0_n$ (see Figure~\ref{fg:1}).  Then whp there exists a subgraph $F_n  \subset G(n, p)$ which is weakly $(G(n, p), H)$-saturated and $F_n$ has $\min\{|E(G(n,p)),|E(F^0_n)|\}$ edges.\\
\label{transfer}
\end{theorem}

%\textbf{Lemma 1 (transferring lemma)} Let $H$ be an arbitrary fixed graph, $p \in (0,1)$. Let $F^0_n$ be a graph on $n$ vertices such that it is weakly$(K_n, H)$-saturated and has the following structure: $F^0_n$ consists of a finite  subgraph $S_0$ and there are $d \ge \min \deg H - 1$ edges to $S_0$ from all other vertices of $K_n$. Then w.h.p. $\exists F  \subset G(n, p)$($F$ is also random) :$F$ is $(G(n, p), H)$-saturated and $F$ has the same number of edges as $F^0_n$.

This theorem implies\\

\begin{corol}
Let $p\in(0,1)$ be constant. For an arbitrary graph $H$ without isolated vertices,  whp 
$$
\mathrm{wsat}(G(n, p), H) = \mathrm{wsat}(K_n, H),
$$ 
if, for every $n\in\mathbb{N}$, there exists a minumum (having $\mathrm{wsat}(K_n, H)$ edges) weakly $(K_n, H)$-saturated graph with the property described in Theorem~\ref{transfer}.\\
\label{corol_transfer}
\end{corol}

\begin{proof}

Indeed, assume that the condition from Theorem~\ref{transfer} is satisfied. Then, it immediately implies that whp $\mathrm{wsat}(G(n, p), H) \le \mathrm{wsat}(K_n, H)$. Assume that, with a probability bounded away from 0, $\mathrm{wsat}(G(n, p), H)$ is strictly less than $\mathrm{wsat}(K_n, H)$. Since whp every pair of vertices of $G(n, p)$ has at least $|V(H)|$ pairwise adjacent common neighbors~\cite{Spencer}, whp $G(n,p)$ is weakly $(K_n, H)$-saturated --- a contradiction.\\ %So, w.h.p. we can use $wsat(G(n, p)$-edged subgraph to restore at first all the edges of $G(n, p)$ and then all the edges of $K_n$. It implies that w.h.p. $wsat(K_n, H) \le wsat(G(n, p), H)$.
\end{proof}

This corollary immediately implies stability for several pattern graphs. Notice that the graph obtained by removing a copy of $K_{n-s+2}$ from $K_n$ is weakly $(K_n, K_s)$-saturated, has the structure described in Theorem~\ref{transfer} and has $\mathrm{wsat}(K_n, K_s)$ edges. Therefore, the result of Kor\'{a}ndi  and Sudakov can be immediately deduced from Corollary~\ref{corol_transfer}.

The constructions of Kronenberg, Martins and Morrison \cite{martins} of weakly $(K_n,K_{t,t})$-saturated and weakly $(K_n,K_{t,t+1})$-saturated graphs with $\mathrm{wsat}(K_n,K_{t,t})$ and $\mathrm{wsat}(K_n,K_{t,t+1})$ edges respectively also have the structure described in Theorem~\ref{transfer}. Therefore, by Corollary~\ref{corol_transfer}, for every $p\in(0,1)$, whp 
$$
\mathrm{wsat}(G(n, p), K_{t, t}) = \mathrm{wsat}(K_n, K_{t,t}),
\quad
\mathrm{wsat}(G(n, p), K_{t, t+1}) = \mathrm{wsat}(K_n, K_{t, t+1}).
$$

The bounds (\ref{general_bipartite_upper}), (\ref{general_bipartite_lower}) imply that Corollary~\ref{corol_transfer} can be also applied to pattern graphs $K_{s,t}$ for all possible values of $s\leq t$. To see this it is sufficient to show that, for every $n\in\mathbb{N}$, there exists a minimum weakly $(K_n,K_{s,t})$-saturated graph with the property described in Theorem~\ref{transfer}. Note that due to (\ref{general_bipartite_upper}) and (\ref{general_bipartite_lower}) $\mathrm{wsat}(K_n,K_{s,t})=(s-1)n+O(1)$, and $s$ is the minimum degree of $K_{s,t}$. Moreover, it is clear that there exists a constant $C=C(s,t)$ such that, for all $n$ large enough, $\mathrm{wsat}(K_n,K_{s,t})=(s-1)n+C$. Indeed, otherwise there exist $C_1<C_2$ and two infinite sequences $\{n^1_i\}_{i\in\mathbb{N}}$ and $\{n^2_i\}_{i\in\mathbb{N}}$ such that
$$
\mathrm{wsat}(K_{n^1_i},K_{s,t})=(s-1)n^1_i+C_1\text{ and }\mathrm{wsat}(K_{n^2_i},K_{s,t})=(s-1)n^2_i+C_2.
$$
Let us choose sufficiently large $i$ and $j$ such that $n_1:=n^1_i<n^2_j=:n_2$ and let $F_1$ be weakly $(K_{n_1},K_{s,t})$-saturated. Then, the graph $F_2$ on $[n_2]$ obtained from $F_1$ by adding $s-1$ edges from each of the vertices from $[n_2]\setminus[n_1]$ to $[n_1]$ is $(K_{n_2},K_{s,t})$-weakly saturated and has the number of edges $(s-1)n^2_i+C_1<(s-1)n^2_i+C_2$ --- a contradiction. From this, the existence of a weakly $(K_n,K_{s,t})$-saturated graph with the desired property is straightforward. Indeed, let $n_0$ be so large that $\mathrm{wsat}(K_{n_0},K_{s,t})=(s-1)n_0+C$ and let $F^0_{n_0}$ be a weakly $(K_{n_0},K_{s,t})$-saturated graph with $(s-1)n_0+C$ edges. For every $n>n_0$, set $S^0_n=[n_0]$ and define $F^0_n$ to be the union of $F^0_{n_0}$ with $s-1$ edges going from each of the vertices from $[n]\setminus S^0_n$ to $S^0_n$. The graph $F^0_n$ is weakly $(K_n,K_{s,t})$-saturated and has the desired number of edges. Therefore by Corollary~\ref{corol_transfer}, we get that, for every constant $p\in(0,1)$ and all $1\leq s\leq t$, whp 
$$
\mathrm{wsat}(G(n, p), K_{s, t}) = \mathrm{wsat}(K_n, K_{s,t}).
$$
This stability result demostates the efficiency of Theorem~\ref{transfer}: it can be used to prove stability even when the exact value of $\mathrm{wsat}(K_n, K_{s,t})$ is unknown.\\

%Corollary~\ref{} can be also applied to show that the weak saturation number is asymptotically stable for complete multipartite graphs. Indeed, in~\cite{} it is shown that 

%\smallskip

%{\it Remark.} % [general bounds for random graph...]n(s-1)+...
%In the same way, for arbitrary $2\leq s<t$, from~(\ref{general_bipartite_upper}) and (\ref{general_bipartite_lower}), we get that there exists $C=C(s,t)$ such that whp
%$$
%|\mathrm{wsat}(G(n,p),K_{s,t})-n(s-1)|<C.
%$$

\medskip

%also described weakly $(K_n, K_{t, t})$ or $(K_n, K_{t, t})$-saturated graphs with $wsat(K_n, K_{t, t})$ or $wsat(K_n, K_{t, t+1})$ edges respectively that have structure described in Lemma 1. Using Corollary 1 one can easily see that w.h.p. $wsat(G(n, p), K_{t, t+1}) = wsat(K_n, K_{t, t})$ and w.h.p. $wsat(G(n, p), K_{t, t}) = wsat(K_n, K_{t, t+1})$ for $p \in (0, 1)$.

However, there are graphs $H$ such that minimum weakly $(K_n,H)$-saturated graphs do not satisfy the condition from Theorem~\ref{transfer}. For such graphs, the same stability result can not be proven using Corollary~\ref{corol_transfer}. 

For example, consider $H$ being the, so called, $t$-barbell graph consisting of two vertex-disjoint $t$-cliques together with a single edge that has an endpoint in each clique. It is clear that $\mathrm{wsat}(K_n,H)\leq {t\choose 2}\frac{n}{t}\sim\frac{tn}{2}$ for $n$ divisible by $t$ since the disjoint union of $n/t$ $t$-cliques is weakly $(K_n,H)$-saturated. Nevertheless, any subgraph of $K_n$ with the property described in Theorem~\ref{transfer} has at least $(t-2)n+O(1)$ vertices. Therefore, it is not minimum possible, and Corollary~\ref{corol_transfer} is not applicable.

Nevertheless, some extra work (it is very technical, so we omit the details) is required to show that the union of disjoint cliques has minimum possible number of edges and that, since whp $G(n,p)$ contains a $K_t$-factor (see~\cite{clique-factors}), we get stability for $t$-barbell graph as well. In some sense, the situation covered by Theorem~\ref{transfer} is less pleasant. Indeed, if a graph on $[n]$ has a bounded maximum degree, then whp $G(n,p)$ contains its isomorphic copy as a spanning subgraph (see~\cite{AlonFuredi}). So, for such weakly saturated graphs, stability is straightforward. In Theorem~\ref{transfer}, we consider an opposite scenario --- $F_n^0$ has vertices with degrees $n-O(1)$ that are not likely to be in $G(n,p)$.  Having that in mind, we conjecture that, for any constant $p\in(0,1)$ and every graph $H$, whp $\mathrm{wsat}(G(n,p),H)=\mathrm{wsat}(K_n,H)$.\\

%Another example of such an $H$ is the $(t,1)$-lollipop graph, i.e. a union of a $t$-clique and a single vertex adjacent to a unique vertex of the clique. For this pattern graph, $\mathrm{wsat}(K_n,H)={t\choose 2}$. However,  \\

Let us now switch to the case $p = o(1)$. 

Kor\'{a}ndi and Sudakov \cite{sudakov} claim that their result can be easily extended to the range $n^{-\varepsilon(s)} \le p \le 1$. It can be  seen that the same is true for Theorem~\ref{transfer} with $\varepsilon$ depending only on $H$. However, for smaller $p$ Theorem~\ref{transfer} or Corollary~\ref{corol_transfer} may not hold.

Kor\'{a}ndi and Sudakov \cite{sudakov} pose the following question: what is the exact probability range where whp $\mathrm{wsat}(G(n, p), K_s) = \mathrm{wsat}(K_n, K_s)$? In 2020, Bidgoli et al. \cite{bidgoli} proved the existence of threshold probability  for  this stability property.  Moreover, they obtained bounds on the threshold:
\begin{itemize}
\item there exists $c$ such that, if $p <c n^{-\frac{2}{s+1}}(\ln n)^{\frac{2}{(s-2)(s+1)}}$, then whp $\mathrm{wsat}(G(n, p), K_s) \neq \mathrm{wsat}(K_n, K_s)$,

\item if $p > n^{-\frac{1}{2s-3}}(\ln n)^2$, then whp $\mathrm{wsat}(G(n, p), K_s) = \mathrm{wsat}(K_n, K_s)$.
\end{itemize}

In this paper, we estimate the threshold probability for the {\it weak $K_{1,t}
$-saturation stability property} $\mathrm{wsat}(G(n, p), K_{1, t}) = \mathrm{wsat}(K_n, K_{1, t})$.

\begin{theorem} Let $t\geq 3$. 
Denote $p(n, t) = n^{-\frac{1}{t-1}}[\ln n]^{-\frac{t-2}{t-1}}$.

\begin{itemize}
\item There exists $c>0$ such that, if $\frac{1}{n^2}\ll p<c p(n, t)$, then whp $\mathrm{wsat}(G(n, p), K_{1, t}) \neq \mathrm{wsat}(K_n, K_{1, t})$.

\item There exists $C>0$ such that, if $p>C p(n, t)$, then whp $\mathrm{wsat}(G(n, p), K_{1, t}) = \mathrm{wsat}(K_n, K_{1, t})$.
%\item If $p\gg n^{-\frac{1}{t-1}}$, then whp $\mathrm{wsat}(G(n, p), K_{1, t}) = \mathrm{wsat}(K_n, K_{1, t})$.
\end{itemize}
\label{threshold_stars}
\end{theorem}

Note that Theorem~\ref{threshold_stars} does not cover the case $t=2$ as well as $p=O(1/n^2)$. But these cases are much easier. Below we consider them separately.

First, if $p < \frac{Q}{n^2}$ for some constant $Q>0$, then whp $G(n,p)$ consists of isolated vertices and isolated edges (it simply follows from Markov's inequality applying to the number of $P_3$ in $G(n,p)$). Therefore, whp there are no copies of $K_{1, t-1}$ for $t \ge 3$ in $G(n, p)$, and there are no weakly $(G(n, p), K_{1, t})$-saturated subgraphs other than the entire graph. So, whp there is stability only if the number of edges of the graph is exactly ${t\choose 2}$. The latter property holds with probability ${{n\choose 2}\choose {t\choose 2}}p^{{t\choose 2}}(1-p)^{{n\choose 2}-{t\choose 2}}$. Therefore, it tends to 0 when $p\ll\frac{1}{n^2}$ and it is bounded away both from 0 and from 1 when $\frac{q}{n^2}<p<\frac{Q}{n^2}$ for some $0<q<Q$.

%{\it Remark.} The lower and upper bounds for the threshold probability in Theorem~\ref{threshold_stars} differs in $[\log n]^{\frac{t-2}{t-1}}$ factor. We conjecture that our lower bound is sharp, i.e. $n^{-\frac{1}{t-1}} [\ln n]^{-\frac{t-2}{t-1}}$ is the threshold probability for the weak saturation stability property.\\

The case $t=2$ is also trivial. Clearly, for a graph $G$ on $[n]$, 
$$
\mathrm{wsat}(G,K_{1,t})=\mathrm{wsat}(K_n,K_{1,t})={t\choose 2}=1
$$ 
if and only if $G$ has exactly one non-empty connected component. Using the standard first and second moment methods (see, e.g.,~\cite[Chapter 1]{Janson_book}), it can be proven that, for every $\varepsilon>0$,
\begin{itemize}
    \item if $p>(1+\varepsilon)\frac{\ln n}{2n}$, then whp $G(n,p)$ contains a unique non-empty component;
    \item if $\frac{1}{n^2}\ll p<(1-\varepsilon)\frac{\ln n}{2n}$, then whp $G(n,p)$ contains at least two non-empty connected components;
    \item if $\frac{q}{n^2}<p<\frac{Q}{n^2}$ for some $0<q<Q$, then whp all edges in $G(n,p)$ are disjoint and, arguing as above, we get that stability happens only if the graph contains ${t \choose 2} = 1$ edges;
    \item if $p\ll\frac{1}{n^2}$, then whp $G(n,p)$ is empty.
\end{itemize}
Therefore, 
\begin{enumerate}
    \item if $p>(1+\varepsilon)\frac{\ln n}{2n}$, then whp $\mathrm{wsat}(G(n,p),K_{1,t})=\mathrm{wsat}(K_n,K_{1,t})$;
    \item if $\frac{1}{n^2}\ll p<(1-\varepsilon)\frac{\ln n}{2n}$, then whp $\mathrm{wsat}(G(n,p),K_{1,t})\neq\mathrm{wsat}(K_n,K_{1,t})$;
    \item if $\frac{q}{n^2}<p<\frac{Q}{n^2}$ for some $0<q<Q$, then 
\begin{multline*}
    {\sf P}\biggl[\mathrm{wsat}(G(n,p),K_{1,t})=\mathrm{wsat}(K_n,K_{1,t})\biggr]=\\
    {\sf P}(G(n,p)\text{ contains exactly 1 edge})+o(1)={n\choose 2}p(1-p)^{{n\choose 2}-1}+o(1)
\end{multline*}
    is bounded away both from 0 and 1,
    \item if $p\ll\frac{1}{n^2}$, then whp $\mathrm{wsat}(G(n,p),K_{1,t})=0\neq\mathrm{wsat}(K_n,K_{1,t})$.
\end{enumerate}

The structure of the paper is the following. In Section 2, we prove Theorem~\ref{transfer}. In Section 3, we prove Theorem~\ref{threshold_stars}.

%respectively. Finally, in Section 5, we make a conclusion and discuss some open questions.

\section{Proof of Theorem~\ref{transfer}}

%Let us now provide proof for lemma \ref{transfer}.

%Let $F$ be some weakly $(G_n, H)$-saturated subgraph. At first, let us prove the equivalence of the following conditions:

%\begin{enumerate}
%    \item There exists $C$ such that $F$ contains a set of vertices $S$ of size at most $C$ such that every vertex of $G_n\setminus S$ is adjacent to  at least $d$ vertices of $S$.
%    \item There exists $D$ such that $F^0_n$ contains a set of vertices $S$ of size at most $D$ such that every vertex of $G_n\setminus S$ is adjacent to  \textbf{exactly} $d$ vertices of $S$.
%\end{enumerate}

We denote $d = \delta(H)-1$, $r = |V(H)|$.  First, for convenience, let us state an equivalent modification of Theorem \ref{transfer} which at first sight seems to be weaker. However, it is equivalent, and it is easier to prove this version. In Section~\ref{sc:equivalence_proof}, we prove the equivalence of the two statements, and then prove the modified version.% is equivalent to Lemma \ref{transfer}.

\begin{lemma}
Let $H$ be an arbitrary fixed graph without isolated vertices, $p \in (0,1)$, $C\geq d$ be constants. For every $n\in\mathbb{N}$, let $F^1_n$ be a weakly $(K_n, H)$-saturated graph containing a set of vertices $S^1_n$ of size at most $C$ such that every vertex of $[n]\setminus S^1_n$ is adjacent to  \textbf{exactly} $d$ vertices of $S^1_n$ and there are no edges between vertices of $[n]\setminus S^1_n$. Then whp there exists a subgraph $F_n  \subset G(n, p)$ which is weakly $(G(n, p), H)$-saturated and $F_n$ has the same number of edges as $F^1_n$.
\label{transfer_1}
\end{lemma}

\subsection{Proof of equivalence}
\label{sc:equivalence_proof}

Clearly, Theorem \ref{transfer} implies Lemma \ref{transfer_1}. Let us prove that the opposite implication is also true.

%Now let us provide the proof. 

%$2 \implies 1$:  obvious.

We first notice that $\mathrm{wsat}(K_n, H) \le {r \choose 2} + d (n-r)$. Indeed, we can construct a weakly $(G_n, H)$-saturated subgraph with at most so many edges in the following way. Let $F_0$ be a weakly $(K_r, H)$-saturated subgraph on $[r]$,  the first $r$ vertices of $K_n$. The desired graph $F_1$  has all the same edges as $F_0$ on $[r]$ and from every other vertex there are exactly $d$ edges to $F_0$. %Notice that $F_1$ is weakly $(G_n, H)$-saturated. Indeed, at first we can restore all edges of $G_n|_{[r]}$. \\Suppose that all edges of $G_n|_{[i]}$ are restored. There are at least $d$ edges from vertex $i+1$ to $G_n|_{[i]}$. Then, for every $u \in G_n|_{[i]}$ there is a copy of $H' \in G_n|_{[i]} \setminus \{u\}$ (let the isomorphism between $H'$ and this copy be $\varphi$) such that for every $w \in H'$ $\ (w, w_1) \in E(H) \implies (\varphi(w), i+1) \in E(F_1)$ and $(\varphi(w), u) \in E(G_n|_{[i]})$ as $G_n|_{[i]}$ is a clique. So, there is a copy of $H''$ where $w_1$ and $w_2$ are mapped to  $i+1$ and $u$ respectively and we can restore the edge $(i+1, u)$. So, as $u$ is arbitrary, we can restore all the edges of $G_n|_{[i+1]}$. 
The existence of a bootstrap percolation process that starts on $F_1$ and finishes on $K_n$ is straightforward: first restore all edges of $K_r$, then restore all edges going to $[r]$ and finally restore all the remaining edges. 

%By the principle of mathematical induction we can restore all the edges of $G_n$. 

We assume that Lemma~\ref{transfer_1} holds. The constructed graph $F_1$ satisfies the conditions of this lemma. Then, whp  there exists a weakly $(G(n, p), H)$ -saturated subgraph $F'_n$ such that 
$$
|E(F'_n)| = |E(F_1)| \leq {r \choose 2} + d (n-r) = dn + \left({r \choose 2} - dr\right).
$$

Let $F^0_n$ be a graph that satisfies the condition of Theorem \ref{transfer}.

If $|E(F^0_n)| \ge |E(F'_n)|$, then a subgraph of $G(n,p)$ obtaining by adding $\min\left\{|E(F^0_n)|,|E(G(n,p))|\right\} - |E(F'_n)|$ edges to $F'_n$ is weakly $(G(n, p), H)$-saturated and $F_n$ has $\min\left\{|E(F^0_n)|,|E(G(n,p))|\right\}$ edges.

If $|E(F^0_n)| < |E(F'_n)|$, then there are at most $a$ vertices with degree more than $d$ in $F^0_n$ outside $S^0_n$. We can add them to $S^0_n$, and then the size of this set will be bounded from above by $C+a$. So, $F^0_n$ satisfies the condition of Lemma \ref{transfer_1} with $C:=C+a$. Therefore, whp there exists a weakly $(G(n, p), H)$-saturated subgraph with $|E(F^0_n)|$ edges. $\Box$\\

Below, we prove Lemma~\ref{transfer_1}. Let us outline the proof. In Section~\ref{properties}, we describe sufficient properties of a spanning subgraph $G$ of $K_n$ that allow to find a weakly $(G,H)$-saturated subgraph with the same number of edges as in a weakly $(K_n,H)$-saturated subgraph. In Section~\ref{properties_are_sufficient}, we prove that these properties are indeed sufficient for the described transference property. In Section~\ref{random_proof}, we prove that whp $G(n,p)$ has the described properties and, thus, finish the proof of Lemma~\ref{transfer_1}.\\

\smallskip

Let us now switch to the proof of Lemma~\ref{transfer_1}. Assume that the requirements of the lemma hold. We also suggest that $n$ is even to avoid overloading with floor and ceiling functions notations. This does not affect the proof anyhow.

Let us denote the vertices of $H$ as $w_1, \ldots, w_{r}$ where $w_1$ is a vertex with degree $d+1$ and $w_2, \ldots, w_{d+2}$ are its neighbours. Let $H'$ be obtained from $H$ by deleting the vertices $w_1, w_2$ and let $H''$  be obtained from $H$ by deleting the edge $\{w_1, w_2\}$ (but preserving the vertices $w_1, w_2$).

In what follows, for a graph $F$ and its vertex $v$, we denote by $N_F(v)$ the set of all neighbours of $v$ in $F$.

\subsection{Sufficient properties}
\label{properties}

We want to find in $G(n,p)$ a subgraph having similar structure to a weakly saturated subgraph in $K_n$. However, it is not immediate since whp all vertices in $G(n,p)$ have degrees $np(1+o(1))$ which is far away from $n-O(1)$. Nevertheless, we can find a clique $K$ in $G(n,p)$ of size $\Theta(\ln n)$, and first reconstruct the edges of the clique. For that, we fix in $K$ a weakly saturated spanning subgraph with the minimum possible number of edges and the desired structure. In other words, we choose a subset $S\subset K$ playing the role of $S^0_{|K|}$. After reconstructing the edges of $K$ we might hope that it is sufficient to use $d$ edges of $G(n,p)$ at every vertex outside $K$ to reconstruct all the other edges of $G(n,p)$. The properties that allow to do this are decribed below.\\

%Suppose that the prerequisites of \ref{transfer} are satisfied. At first, let us formulate a property of a non-random graph $G$ on $[n]$ which implies that this graph has a weakly $(G, H)$-saturated subgraph with at most $|E(F_n)|$ edges. After that we will prove that $G(n, p)$ whp satisfies this property.\\
We start from distinguishing several subsets of $[n]$ that we use to describe the properties. Everywhere below, by $G|_V$ we denote the induced subgraph of $G$ on the vertex set $V$ (where $V \subset V(G)$).\\

Let $c>0$. Let $G$ be a graph on the vertex set $[n]$. Let
\begin{itemize}
\item $G_1=G|_{[n/2]}$, $G_2=G|_{[n]\setminus[n/2]}$; 

\item $K\subset V_1:=V(G_1)$ be a set of size $k\geq c\ln n$, 

$S$ be a subset of $K$ of size $|S^1_k|$, 

$D$ be a subset of $S$ of size $d$ (clearly, the requirements in Lemma~\ref{transfer_1} imply that $|S^1_k|\geq d$);

\item $Z$ be the set of all common neighbours of $D$ in $V(G_2)$;

\item $Z=Z_1\sqcup Z_2\sqcup Z_3$ be a partition of $Z$. %into subsets of almost equal sizes, i.e. $||Z_i|-|Z_j||\leq 1$.%, $||Z^1|-|Z^2||\leq 1$. \textcolor{blue}{//second partition is not equal, it will be random}

\item $R$ be a set of $r$ vertices from $K \setminus S$ and $T\subset Z$ be the set of all common neighbors from $Z$ of vertices from $R$. %?? ??????. ????? ????, ???????? ??? ???? ?? ???????:

%\textcolor{blue}{For every $u, v \in T$ there exists a $H^1 \cong H' \in K \setminus S$ such that $(u, v, H^1)$ is $H$-completable}.

\end{itemize}

%For $v\in Z\setminus T$, denote by $U_v$ the set of all neighbors of $v$ in $T$.% having at least $d$ common neighbors with $v$ in $K\setminus S$.

\begin{figure}
\begin{centering}
\begin{tikzpicture}[main/.style = {node distance={15mm}, draw, circle}] \small

\filldraw[color=red!90, fill=red!05] (0.5, 0.5) circle (1.2);

\node[] () [above of=3] {$H$}; 
\node[red] () [above of=1] {$H'$}; 
\node[main] (1) {}; 
\node[main] (2) [above right of=1] {}; 
\node[main] (3) [below left of=1]{$w_1$}; 
\node[main] (4) [below right of=1]{$w_2$}; 
\node[main] (5) at (1.065, 0) {}; 

\draw[-] (1) -- (2);
\draw[-] (1) -- (3);
\draw[-] (1) -- (4);
\draw[-] (3) -- (4);
\draw[-] (2) -- (5);
\draw[-] (4) -- (5);

\filldraw[color=gray!90, fill=gray!0] (7, 0) circle [x radius=3.5cm, y radius=2.5cm];

\filldraw[color=red!90, fill=red!05] (6.5, 0.5) circle (1.2);

\filldraw[color=violet!50, fill=violet!25] (6.55, 0) circle [x radius=0.9cm, y radius=0.35cm];
\node[violet] at (7.75, 0)  () {$\tilde V$};

\node[] at (9, 0)  () {$G$};

\node[main] at (6, 0)  (1) {}; 
\node[red] () [above  of=1] {$\tilde H$}; 
\node[main] (2) [above right of=1] {}; 
\node[main] (3) [below left of=1]{$v$}; 
\node[main] (4) [below right of=1]{$u$}; 
\node[main] (5) at (7.065, 0) {}; 
\draw[-] (1) -- (2);
\draw[-] (1) -- (3);
\draw[-] (1) -- (4);
\draw[-] (2) -- (5);

\draw[dash pattern=on 4pt off 2pt, line width=3, blue] (3) -- (4);
\draw[dash pattern=on 2pt off 2pt] (5) -- (1);
\draw[-] (4) -- (5);

\node[main, dash pattern=on 2pt off 2pt] (a) at (8.5, 0.575) {}; 
\node[main, dash pattern=on 2pt off 2pt] (b) at (8.5, 1.5) {}; 

\draw[dash pattern=on 2pt off 2pt] (5) -- (a);
%\draw[dash pattern=on 2pt off 2pt] (5) -- (b);
\draw[dash pattern=on 2pt off 2pt] (2) -- (a);
\draw[dash pattern=on 2pt off 2pt] (2) -- (b);
%\captionof{figure}{$H$-completable tuples $(v, u, \tilde H)$ and $(v, u, \tilde H, \tilde V)$. Black edges are edges mapped to edges of $H$, dashed edges are other edges of $G$, the blue dashed edge can be immediately added to $G$ in a bootstrap percollation process.}
\end{tikzpicture} 

\caption{$H$-completable tuples $(v, u, \tilde H)$ and $(v, u, \tilde H, \tilde V)$. Black edges are edges mapped to edges of $H$, dashed edges are other edges of $G$, the blue dashed edge can be immediately reconstructed in a bootstrap percollation process.}
\label{fg:2}

%\textcolor{blue}{In Figure 2, an example of an $H$-completable tuple can be seen. Black edges are edges mapped to edges of $H$, dashed edges are other edges of $G$, the blue dashed edge can be immediately added to $G$ in a bootstrap percollation process.}
\end{centering}
\end{figure}
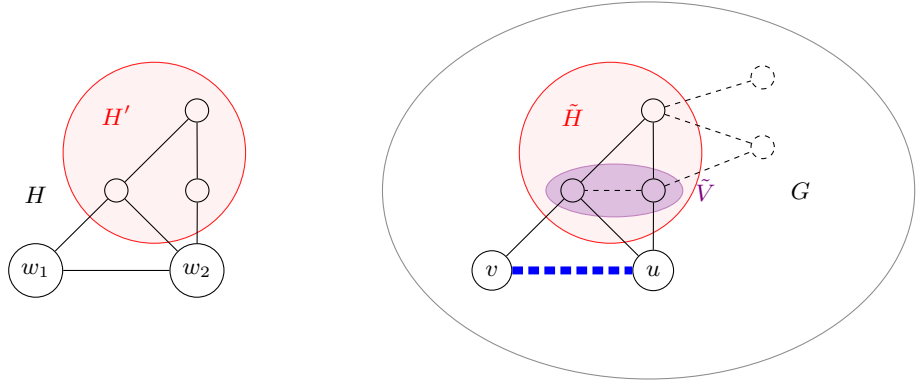

For vertices $v,u$ of $G$ and a subgraph $\tilde H\cong H'$ of $G$, we call the tuple $(v,u,\tilde H)$ {\it $H$-completable in $G$} (see Figure~\ref{fg:2}), if there exists an embedding $f$ (we call it {\it $(v,u,\tilde H)$-embedding}) from $H$ to $G|_{V(\tilde H)\sqcup\{v,u\}}$ such that $f(w_1)=v$, $f(w_2)=u$ and $f$ maps $H'$ to $\tilde H$, i.e. the graph with the set of vertices $V(\tilde H)$ and the set of edges $\{\{f(x),f(y)\},\{x,y\}\in E(H')\}$ equals $\tilde H$. In plain words, it means that we are able to immediately reconstruct the edge $\{u, v\}$. % to $G$, when $G$ appears during a bootstrap percolation process.}
For a vertex $v$ of $G$ and a subgraph $\tilde H\cong H'$ of $G$, we call the pair $(v, \tilde H)$ {\it $H$-completable in $G$}, if there exists an embedding $f$ (we call it {\it $(v,\tilde H)$-embedding}) from $H|_{V(H)\setminus\{w_2\}}$ to $G|_{V(\tilde H)\sqcup\{v\}}$ such that $f(w_1)=v$ and $f$ maps $H'$ to $\tilde H$. Intuitively, it means that nothing prevents us from adding edges from $v$ to some other vertices $w$ (such that $(v, w, \tilde H)$ is  $H$-completable in $G$). Although we do not know whether such $w$ exists, it is not 'the fault' of $v$. Let $(v,u,\tilde H)$ be $H$-completable in $G$ and $\tilde V \subset V(\tilde H)$ have exactly $d$ vertices. We call the tuple $(v,u,\tilde H, \tilde V)$ {\it $H$-completable in $G$}, if there exists a $(v,u,\tilde H)$-embedding $f$ such that $f$ maps $N_H(w_1) \setminus \{w_2\}$  onto $\tilde V$. In Figure~\ref{fg:2}, $\tilde V$ is shown in violet color. Similarly, for an $H$-completable $(v,\tilde H)$ and $\tilde V \subset V(\tilde H)$ of size $d$, we call the tuple $(v,\tilde H, \tilde V)$ {\it $H$-completable in $G$}, if there exists an $(v,\tilde H)$-embedding $f$ such that $f$ maps $N_H(w_1) \setminus \{w_2\}$  onto $\tilde V$. \\

% vertex $v$ of $G$, subgraph $\tilde H\cong H'$ of $G$ and set of vertices $\tilde V \subset V(\tilde H)$ of size $d$, we call the tuple $(v,u,\tilde H, \tilde V)$ {\it $H$-completable in $G$}, if there exists an embedding $f$ from $H|_{V(H)\setminus w_2}$ to $G|_{V(\tilde H)\sqcup\{v\}}$ such that $f(w_1)=v$,  $f$ maps $H'$ to $\tilde H$, i.e. the graph with the set of vertices $V(\tilde H)$ and the set of edges $\{\{f(x),f(y)\},\{x,y\}\in E(H')\})$ equals $\tilde H$, and $f$ maps $N_H(w_1) \setminus \{w_2\}$  onto $\tilde V$.}

%\textcolor{blue}{For $v \in Z\setminus T$, denote by $U_v$ the set of all neighbors of $v$ in $T$ such that there exists $H_v\cong H'$ inside $K \setminus S$ with $(v,u,H_v)$ being $H$-completable.} 
%
%\textcolor{blue}{Now, let us define a tuple $(G,c,K,S,D,Z_1,Z_2,Z_3,T)$. Here, $G$ is a graph, $c$ is a constant, $K,S,D,Z_1,Z_2,Z_3,T$ are subsets of $V(G)$. We also assume that $K$ and $Z = Z_1 \sqcup Z_2 \sqcup Z_3$ are disjoint, $D \subset S \subset K$, $T \subset Z$ (but there is no other relation between $T$ and $Z_i, i \in \{1, 2, 3\}$.  Intuitively, we want to create in a random graph a structure similar to the structure of the weakly-saturated subgraph in the full graph. Here, $S$ will play the part of $S^0_n$ in the construction of the full graph. However, the rest of the construction will be more complicated. }  

Let us now desribe the desired properties. We define two properties of $G$. The first one is used to reconstruct all the edges but those between $S\setminus D$ and $[n]\setminus(K\sqcup Z)$. All the other edges are reconstructed using the second property.

For a vertex $v$ of $G$, a subgraph $\tilde H\cong H'$ of $G$ and a set of vertices $\tilde V \subset V(\tilde H)$ of size $d$, denote by $U_v(\tilde H, \tilde V)$ the set of all neighbours $u$ of $v$ in $Z_2$ such that the tuple $(v,u,\tilde H, \tilde V)$ is $H$-completable. 
\\

%For a vertex $v$ of $G$, subgraph $\tilde H\cong H'$ of $G$ denote by $U_v(\tilde H)$ the set of all neighbours $u$ of $v$ in $Z_2$ such that the tuple $(v,u,\tilde H)$ is $H$-completable.\\

Let us say that the tuple $(G;K,S,D,Z_1,Z_2,Z_3,T)$ satisfies the {\it first $H$-saturation property} (or, simply, {\it the first property}), if (note that, among the members of the tuple, $G$ is a graph, and all the others are sets of vertices; recall also that $Z_1\sqcup Z_2\sqcup Z_3$ is a partition of $Z$)
\begin{enumerate}
    \item \label{K} $K$ induces a clique in $G$;
%    \item \label{size} $|Z|\geq cn$, $|T|\geq n^c$;
%    \item  \label{D} there are at least $c^2 n$ common  neighbors of $D$ in $Z$;
    %\item \label{construction_need} each vertex outside $K\sqcup Z$ has at least $d$ neighbors in $Z$;
    \item \label{Z_T} {for any adjacent (in $G$) pair $v_1,v_2\in Z\setminus T$, there exists a copy $H_{v_1, v_2}$ of $H'$ in $G|_T$ such that $(v_1, v_2, H_{v_1, v_2})$ is $H$-completable in $G$};
    %\item  \label{Z_T&T} for any 
    %$v\in Z\setminus T$, $u\in T$, there exists a copy of $H'$ in $G|_{U_{v}\cap N_G(u)}$;

    %\item \label{H_v} {for every $v\notin K\cup Z$ there exists $H_v \cong H'$ inside $Z_1$ such that $(v, H_v)$ is $H$-completable in $G$; }
    
    \item \label{H_v} for every $v\notin K\sqcup Z$ there exists $H_v \cong H'$ inside $Z_1$ and $V_v \subset V(H_v)$ such that $(v, H_v, V_v)$ is $H$-completable in $G$; 
    
    \item  \label{others} for any pair $v_1,v_2\notin K\sqcup Z$, there exists a copy of $H'$ in $G|_{U_{v_1}(H_{v_1}, V_{v_1})\cap U_{v_2}(H_{v_2}, V_{v_2})}$;
    \item \label{part_between}for any $v\notin K\sqcup Z$, $u\in (K\setminus S)\sqcup Z_2\sqcup Z_3$, there exists a copy of $H'$ in $G|_{U_{v}(H_v,V_v)\cap N_G(u)}$;
    \item \label{Z1} for any $v\notin K\sqcup Z$, $u\in Z_1$, there exists a copy of $H'$ in $G|_{Z_3\cap N_G(v)\cap N_G(u)}$;
    \item \label{CD} for any $u\notin S$, there exists a copy of $H'$ in $G|_{Z \cap N_G(u)}$.\\
   % \item  \label{Z1D} for any $v\in D$, $u\in Z^1$, there exists a copy of $H'$ in $G|_{Z^2\cap N_G(v)\cap N_G(u)}$.

\end{enumerate}

Finally, let $V(G_2)=V^1\sqcup V^2\sqcup V^3$ be a partition. Let us say that the tuple $(G;V^1,V^2,V^3,S,D)$ satisfies the {\it second $H$-saturation property}  (or, simply, {\it the second property}), if

 %????? ??????? ???, ?? ?????? ???? ???? ??? ???? v ????? ???????  ????? ???? ??????? ?????.

\begin{enumerate}
 %   \item \label{HS} there exists $H_D\subset G|_{D \cup V^1}$ and an isomorphism $H'\to H_D$ that maps $N_H(w_1) \setminus \{w_2\}$ into $D$;
    
   % \textcolor{blue}{For $v \in S \setminus D$, denote by $\hat U_v$ the set of all neighbors of $v$ in $V^2$ such that $(v,u,H_S)$ is $H$-completable.}
    
    \item \label{S_D} there exists $H_D\subset G|_{D \cup V^1}$ such that $H_D \cong H'$ and, %an isomorphism $f:H'\to H_D$ that maps $N_H(w_1) \setminus \{w_2\}$ into $D$ such that, 
    for every $v \in S \setminus D$, $u \in [V(G_1) \setminus S] \sqcup V^3$, there is a copy of $H'$ inside $G|_{\hat U_v \cap N_G(u)}$, where
    $\hat U_v$ is the set of all neighbors $z$ of $v$ in $V^2$ such that $(v,z,H_D, D)$ is $H$-completable.
    
%\begin{center}    
%    $\hat U_v=\{y\in N_G(v)\cap V^2:$ there exists an embedding $\tilde f:H\to G|_{V(H_D)\sqcup\{v,y\}}$ with $\tilde f(w_1)=v$, $\tilde f(w_2)=u$ and $\tilde f|_{H'}=f\}$;
%\end{center}

%    \textcolor{blue}{$(v, y, H_D)$ is $H$-completable.}%$H$ can be embedded in $G|_{V(H_D)\sqcup\{v,y\}}$ in a way such that $w_1$ is mapped to $v$ and $w_2$ is mapped to $y$ and $w_3,\ldots,w_{d+2}$ are mapped onto $D$;
    
    \item \label{S_D_G12} for every $v \in S \setminus D$, $u \in V^1 \sqcup V^2 $, there is a copy of $H'$ inside $G|_{V^3 \cap N_G(u) \cap N_G(v)}$.

\end{enumerate}

\subsection{Proof for a non-random graph}
\label{properties_are_sufficient}

Assume that tuples $(G;K,S,D,Z_1,Z_2,Z_3,T)$ and $(G;V^1,V^2,V^3,S,D)$ satisfy the first and the second property respectively. 

Clearly, in $G$, there are at least ${k\choose 2}+d(n-k)>|E(F^1_n)|$ (for $n$ large enough, since $F^1_n$ has $dn+O(1)$ edges due to the requirements in Lemma~\ref{transfer_1}, and $k\geq c\ln n$) edges. Let us prove that there exists a weakly $(G, H)$-saturated graph with $|E(F^1_n)|$ edges for sufficiently large $n$. Clearly, it is sufficient to prove the existence of a weakly $(G, H)$-saturated graph with at most $|E(F^1_n)|$ edges.

%Let us show that (without loss of generality) we may assume that, for all $n$ large enough, $|V(S^1_n)|=C$ is constant. Assume that this is not the case. Then there exist $n_1<n_2$ such that 
%$$
%e_1:=\left|E\left(F^1_{n_1}|_{S_{n_1}^1}\right)\right|<\left|E\left(F^1_{n_2}|_{S_{n_2}^1}\right)\right|=:e_2.
%$$
%Then $\left|E\left(F^1_{n_1}\right)\right|=e_1+(n_1$

%????? ???????? n_1<n_2 ?????, ??? ? S_{n_1}^0 ????? ?????? ??? ? S_{n_2}^0. ????? ? S_{n_1}^0 e_1 ?????, ? ? S_{n_2}^0 - e_2 ?????. ????? wsat(K_{n_1})=e_1+(n_1-v)d, wsat(K_{n_2})=e_2+(n_2-v)d. ???????, ??? ? K_{n_2} ?????????? weakly saturated ??????? ? ??????? ?????? ?????, ??-?? ???? ?????? ? ????????????. ??????? ? K_{n_2} ??????? K_{n_1}, ? ? ???, ? ???? ???????, weakly saturated ??????? ? e_1+(n_1-v)d ???????. ????? ?? ?????? ?? ?????????? n_2-n_1 ?????? ???????? ? ???? ??????? ?? d ?????. ??????????, ??? ????? ? K_{n_1} ?? ???????????, ????? ???? ??????????? ??? ????? ? K_{n_2}. ?? ? ???????????? ???????? ????? e_1+(n_1-v)d+(n_2-n_1)d=e_1+(n_2-v)d<e_2+(n_2-v)d - ????????????.

Without loss of generality, assume that, for every $n$,  $|V(S^1_n)|=C$ (we can add $C-|V(S^1_n)|$ vertices of $[n]\setminus S^1_n$ to $S^1_n$, and then we will preserve conditions of Lemma~\ref{transfer_1}). We can also assume that, for every $n$, $F^1_n$ has the minimum number of edges among all graphs that satisfy conditions of  Lemma~\ref{transfer_1} and have $|V(S^1_n)| = C$.\\

Let us now construct a spanning subgraph $F\subset G$ with {\it at most} $|E(F^1_n)|$ edges. After that, we will prove that this graph is weakly $(G,H)$-saturated.

%whp  all edges of $G$ can be restored via a bootstrap percolation process from $F$. 

Let us first define those edges of $F$ that belong to $K$.  %\textcolor{blue}{Remember that $F^1_k$ is the weakly-saturated subgraph for $K_k$, and $S^1_k$ is the part of $F^1_k$ of size $C$).} 
Let $\varphi$ be a bijection from $K$ to $V(K_k)$ such that $S\subset K$ is mapped onto $S_k^1$. Then, we construct a graph on $K$ (and we let $F|_K$ to be exactly this graph) isomorphic to $F^1_{k}$ such that $\varphi$ is an isomorphism of $F|_K$ and $F^1_k$ (i.e., an edge $\{u, v\}$ belongs to $F|_K$ iff $\{\varphi(u), \varphi(v)\}$ belongs to $F^1_k$).

Second, for every vertex $v\in Z$, keep the edges of $G$ going from $v$ to $D$ (there are $d$ of them). Moreover, for every vertex outside $K\sqcup Z$, keep {\it specific} $d$ edges of $G$ going from $v$ to $Z$ (so many edges exist due to the condition \ref{H_v} of the first property). The choice of edges between $[n]\setminus(K\sqcup Z)$ and $Z$ will be explained later.\\

Clearly, 
\begin{equation}
|E(F)| = |E(F^1_n)| + \left|E\left(F^1_k|_{S^1_k}\right)\right| - \left|E\left(F^1_n|_{S^1_n}\right)\right|.
\label{edges_in_F}
\end{equation} 
Let us prove that $\left|E\left(F^1_k|_{S^1_k}\right)\right| = \left|E\left(F^1_n|_{S^1_n}\right)\right|$ for $n$ large enough. It is enough to prove that there exists $N \in \mathbb{N}$ such that, for every $n_1, n_2 \ge N$, $\left|E\left(F_{n_1}^1|_{S^1_{n_1}}\right)\right| = \left|E\left(F^1_{n_2}|_{S^1_{n_2}}\right)\right|$. Assume the contrary: for every $N \in \mathbb{N}$, there exist $n_1>n_2 \ge N$ such that $\left|E\left(F_{n_1}^1|_{S^1_{n_1}}\right)\right| > \left|E\left(F^1_{n_2}|_{S^1_{n_2}}\right)\right|$. Let $N$ be large enough. Since $|S^1_{n_1}|=|S^1_{n_2}|=C$, we get that $|E(F_{n_1}^1)|>|E(F_{n_2}^1)|+d(n_1-n_2)$. Since $F^1_{n_2}$ is weakly $(K_{n_2},H)$-saturated, we get that a graph on $[n_1]$ obtained from $F^1_{n_2}$ by adding $d$ edges from each vertex of $[n_1]\setminus[n_2]$ to $F^1_{n_2}$ is weakly $(K_{n_1},H)$-saturated. This contradicts with the minimality of the number of edges in $F^1_{n_1}$.

%First suppose that there exist $N \in \mathbb{N} > |H|, n_1 > n_2 \ge N$ $E(S^1_{n_1}) > E(S^1_{n_2})$. Suppose that in  $F^1_{n_1}$ and $F^1_{n_2}$ the subgraphs $S^1_{n_1}$ and $S^1_{n_2}$ are constructed on $|V(S^1_{n_1})|$ and $|V(S^1_{n_2})|$ respectively, and vertices with larger numbers only have $d$ edges to $S^1_{n_1}$ or $S^1_{n_2}$. Let $\hat F^1_{n_1}$ be a subgraph such that $\hat F^1_{n_1}|_{[n_2]} = F^1_{n_2}$ and from every vertice of $\{n_2+1, \ldots n_1\}$ there are som $d$ edges to $[n_2]$. Notice that $\hat F_{n_1}$ is weakly $(K_{n_1}, H)$-saturated, as we can first restore the edges of $G_{[n_2]} \cong K_{n_2}$ from $F_{n_2}$ and then restore $G_{[n_1]}$ from $G_{[n_2]}$ and $d$ edges from every other vertice. $|E(\hat F_{n_1})| = E(S^1_{n_2}) + d (n - V(S^1_{n_2})) = E(S^1_{n_2}) + d (n - C) < E(S^1_{n_1}) + d (n - C) = |E(F^1_{n_1})|$ which contradicts the fact that $F^1_{n_1}$ has the minimum number of edges among suitable subgraphs.

%But then for every  $N \in \mathbb{N}$ and for every $n_1 > n_2 \ge N$ there is the inequality $E(S^1_{n_1}) \le  E(S^1_{n_2})$ and also there exist $n_1 > n_2 \ge N$ such that $E(S^1_{n_1}) <  E(S^1_{n_2})$. It is impossible, because  $E(S^1_{n})$ can not become below zero. 

From~(\ref{edges_in_F}), we get that $|E(F)| = |E(F^1_n)|$ for large $n$.\\

Now let us show that $F$ is weakly $(G,H)$-saturated and, on the way, specify the edges from $[n]\setminus(K\sqcup Z)$ to $Z$.\\

{\it We first sequentially add  the following bunch of edges to $F$: edges inside $K$,  edges from $K$ to $Z$,  edges inside $T$, edges between $Z$ and $T$, edges inside $Z$.% (here, a few more complicated steps are needed).
}

\begin{itemize}

\item[1)] Here, we restore the edges of $G$ that are inside $K$. This is straightforward since $K$ is a clique (by the condition \ref{K} of the first property), $F|_K \cong F_{k}^1$ and there exists a bootstrap percolation process that starts on $F_k^1$ and finishes on $K_k$. Let $F_1=F\cup G|_K$.

\item[2)] Let us restore the edges of $G$ between $K$ and $Z$. Edges between $Z$ and $D$ are already in $F_1$. Consider $u \in K \setminus D, v\in Z$. Let $K'$ be a set of $r - d - 2$ vertices of $K \setminus [D \sqcup \{u\}]$. Then, for any graph $\tilde H\cong H'$ on the vertex set $K'\sqcup D$ such that $N_H(w_1)$ is mapped onto $D$ % and $w_1$ is mapped onto $v$ %and $N_{G|_{\tilde H\cup \{v\}}}(v)=D$
(such a mapping exists since $K$ induces a clique in $G$), %we can just map $V(H')$  onto $K' \sqcup D$ in any arbitrary way such that $N_H(v)$ is mapped onto $D$ and then take $\tilde H$ as the image of $H'$
 the tuple $(v,u,\tilde H)$ is $H$-completable in $(V(F_1),E(F_1)\sqcup\{u,v\})$ since $v$ is adjacent to every vertex from $D$ in $G$. So, we can restore $\{u, v\}$. Let $F_2$ be obtained from $F_1$ by adding all edges from $G$ between $K$ and $Z$.

\item[3)] Let us switch to edges that are entirely in $T$. Consider $u, v\in T$. The edges inside $R$ and between $\{u, v\}$ and $R$ are already in $F_2$ (since $T\subset Z$, $R\subset K$, and the edges inside $K$ and between $K$ and $Z$ are already restored). Let $\tilde H\cong H'$ be inside $R$ (recall that $R$ is a clique). Then $(v, u, \tilde H)$ is $H$-completable in $(V(F_2),E(F_2)\sqcup\{u,v\})$ and, therefore, we are able to restore $\{u, v\}$. We get $F_3=F_2\cup G|_T$.

\item[4)] Let us restore edges between $Z \setminus T$ and $T$. %For every $v \in Z \setminus T$, we need to restore edges between $v$ and $U_v$ (remember that $U_v = N_G(v) \cap T$). 
Consider $v \in Z \setminus T$, $u \in T$. %U_v = T \cap N_G(v)$. %Notice that from both $v$ and $u$ there are edges to all vertices in $D$ and these edges are in $F_3$.\
Since $u$ is adjacent to all vertices from $R\cup D$, $v$ is adjacent to all vertices in $D$, $|R|=r$ and $|D| = d$, we get that there exists $\tilde H\cong H'$ in $F_3|_{R\cup D}$ such that $(v, u, \tilde H)$ is $H$-completable in $(V(F_3),E(F_3)\sqcup\{u,v\})$. Let $F_4$ be obtained from $F_3$ by adding all edges from $G$ between $Z \setminus T$ and $T$. %every $v\in Z\setminus T$ and $U_v$.

\item[5)] Let us restore the remaining edges inside $Z$. Consider adjacent (in $G$) %First, consider
$v_1,v_2 \in Z \setminus T$. By the condition \ref{Z_T} of the first property, there is $H_{v_1, v_2}\cong H'$ inside $G|_T = F_4|_T$ such that $(v_1, v_2, H_{v_1, v_2})$ is $H$-completable. Since all edges from $G|_T$ and between $\{v_1,v_2\}$ and $T$ are already in $F_4$, we get that $(v_1, v_2,H_{v_1, v_2})$ is $H$-completable in $(V(F_4),E(F_4)\sqcup\{v_1,v_2\})$ and we can restore $\{v_1, v_2\}$. We get $F_5=F_4\cup G|_Z$.\\

{\it Next we restore edges that are entirely outside $K\sqcup Z$.}

\item[6)] Let $v\notin K\sqcup Z$. Notice that {\it we have to specify} $d$ edges going from $v$ to $Z$ in $F$. Let us do that. By the condition \ref{H_v} of the first property, there exists a copy $H_v  \subset G|_{Z_1}$ of $H'$ such that $(v,H_v)$ is $H$-completable. We {\it specify} $V_v$. %the set of all $d$ edges from $v$ to $H_v$ in $G$ corresponding to the edges between $w_1$ and $H'$ in $H$. Let $V_v \subset V(H_v)$ be the set of the second ends of these edges.

Let us now restore edges between $v$ and $U_v(H_v, V_v)$. The edges inside $H_v \subset Z$ and between $v$ and $V_v$ are already in $F_5$, so by the definition of $U_v(H_v, V_v)$, we can restore all edges between $v$ and $U_v(H_v, V_v)$.  Let $F_6$ be obtained from $F_5$ by adding all edges between every $v\notin K \sqcup Z$ and $U_v(H_v, V_v)$.

%By first property \label{others} for any pair $v_1,v_2 \notin K\sqcup Z$, there exists a copy of $H'$ in $G|_{U_{v_1}\cap U_{v_2}}$. 
%It implies that for every $v \notin K \sqcup Z$ the set  $\tilde{U_v}$ is not empty and therefore there exists a copy $H^1_v  \subset Z_1$ of $H'$ from the definition of $\tilde{U_v}$ such that edges from $v$ to $H^1_v$ go the same way as in $H$. 
%Bu first property \ref{H_v} there exists exists a copy $H^1_v  \subset Z_1$ of $H'$ such that $(v, \tilde H_v) $ is $H$-completable.

%Consider v $\notin K \sqcup Z$. Let $E_v$ be the set of edges from $v$ to $\tilde H_v$ corresponding to the edges between $w_1$ and $H'$ in $H$. Notice that $|E_v| = d$. Remember that while constructing $F$ we wanted to draw some $d$ edges from vertices of $G \setminus K \setminus Z$  to $Z$, but the choice of edges was not given. Let these edges be $E_v$ for every $v \notin K \cup Z$.

%Now notice that edges inside $\tilde H_v$ are already restored and there are $d$ correct edges from $v$ to $\tilde H_v$ in $F$. So, by definition of $U_v$ we can restore all edges between $v$ and $U_v(\tilde H_v)$.

\item[7)] Here, we consider edges that have both vertices outside $K\sqcup Z$. Consider $v_1, v_2 \notin K \sqcup Z$. Edges from $v_1$ and $v_2$ to $U_{v_1}(H_{v_1}, V_{v_1}) \cap U_{v_2}(H_{v_2},V_{v_2})$ and edges inside  $U_{v_1}(H_{v_1}, V_{v_1}) \cap U_{v_2}(H_{v_2},V_{v_2})$ are already in $F_6$. By the condition \ref{others} of the first property, there is a copy of $H'$ inside $G|_{U_{v_1}(H_{v_1},V_{v_1}) \cap U_{v_2}(H_{v_2},V_{v_2})}=F_6|_{U_{v_1}(H_{v_1}, V_{v_1}) \cap U_{v_2}(H_{v_2}, V_{v_2})}$, so the edge between $v_1$ and $v_2$ can be restored. Let $F_7=F_6\cup G|_{[n]\setminus[K\sqcup Z]}$.\\

{\it It remains to restore only  edges between $K\sqcup Z$ and $[n]\setminus(K\sqcup Z)$.}

\item [8)] Let us restore all edges between $(K\sqcup Z)\setminus S$ and $[n]\setminus(K\sqcup Z)$. Let $v\notin K\sqcup Z$.

First, let $u \in (K \setminus S) \sqcup Z_2 \sqcup Z_3$. The edges inside $U_v(H_v, V_v) \cap N_G(u) \subset Z$, edges from $u$ to $N_G(u) \cap Z_2$ and edges from $v$ to $U_v(H_v, V_v)$ are already in $F_7$. By the condition \ref{part_between} of the first property, there is a copy of $H'$ inside $U_v(H_v, V_v) \cap N_G(u)$. So, we can restore the edge between $u$ and $v$. 

Second, let $u \in Z_1$. The edges from $v$ to $Z_3$ are just restored. The edges inside $Z_3$ and the edges between $u$ and $Z_3$ are already in $F_7$. By the condition \ref{Z1} of the first property, there is a copy of $H'$ inside $Z_3 \cap N_G(u) \cap N_G(v)$, so we can restore $\{u, v\}$.

$F_8$ is obtained from $F_7$ by adding all edges of $G$ between $(K\sqcup Z)\setminus S$ and $[n]\setminus(K\sqcup Z)$.\\ 

{\it It remains to restore only edges between $S$ and $[n]\setminus(K\sqcup Z)$.}

\item[9)] Here, we restore edges between $D$ and $[n]\setminus(K\sqcup Z)$. Let $v\in D$, $u \in [n]\setminus(K\sqcup Z)$. Then edges between $u$ and $Z$, edges between $v$ and $Z$ and edges inside $Z$ are in $F_8$, and $Z \subset N_G(v)$. By the condition \ref{CD} of the first property, there is a copy of $H'$ inside $Z \cap N_G(u)$. So, we can restore $\{u, v\}$. Let $F_9$ be obtained from $F_8$ by adding all edges of $G$ between $D$ and $[n]\setminus(K\sqcup Z)$.

%Let $u \in Z^1$. All edges between $u$ and $Z^2$, $v$ and $Z^2$ and $Z^2$ are in $F_8$. By first property \ref{Z1D} there is a copy of $H'$ inside $Z^2 \cap N_G(u) \cap N_G(v)$, so we can restore the edge $(u, v)$.

\item[10)] It remains to restore edges between $S\setminus D$ and $[n]\setminus(K\sqcup Z)$. Consider $v \in S \setminus D$. By the definition of $\hat U_v$ (given in the condition \ref{S_D} of the second property), the edges from $v$ to $\hat U_v$ can be restored immediately, as edges inside $H_D$, edges between $\hat U_v$ and $H_D$ and edges from $v$ to $D$ are in $F_9$. For $u \in (V(G_1) \setminus S) \sqcup V^3$, by the condition \ref{S_D} of the second property, there is a copy of $H'$ inside $\hat U_v \cap N_G(u)$. Edges inside $\hat U_v \cap N_G(u) \subset V^2$ and edges between $\{u, v\}$ and $\hat U_v \cap N_G(u) \subset V^2$ are already restored, so $\{u, v\}$ can be restored. Finally, consider $u \in V^1\sqcup V^2$. The edges between $S$ and $V^3$ have just been  restored. By the condition \ref{S_D_G12} of the second property, there is a copy of $H'$ inside $V^3 \cap N_G(u) \cap N_G(v)$, so we can restore $\{u, v\}$ as well.

\end{itemize}

%Property \ref{T} implies that there is a copy of $H'$ inside $K \setminus S$ (already restored) which has ''right'' edges to $u$ and $v$ (also already restored), so the edge $(u, v)$ can be restored (if needed).

%$H_v$ and $U_v$ are as defined in properties \ref{H_in_K} and \ref{Z_T}. Notice that edges inside $H_v$ and between $H_v$ and $U_v$ are already restored, so by definition of $U_v$ all edges from $v$ to $U_v$ can be immediately restored.\\

\subsection{Random graph has the properties}
\label{random_proof}

Let us first recall some results on the distribution of small subgraphs in the binomial random graph.\\

Given a graph $Y$, it is well known that the number of subgraphs isomorphic to $Y$ in $G(n,p)$ is well-concentrated around its expectation. In particular, Janson's inequality implies that (see, e.g.,~\cite[Theorem 2.14]{Janson_book}) the probability that $G(n,p)$ does not contain an isomorphic copy of $K_{\ell}$ ($\ell$ is a positive integer constant) is at most $e^{-\Omega(n^2)}$. By the union bound, we get

\begin{claim}
Let $\varepsilon>0$. Whp, for any subset $A\subset[n]$ such that $|A|\geq\varepsilon n$, there exists a copy of $K_{\ell}$ in $G(n,p)|_A$.
\label{sequence}
\end{claim}

Since $K_{\ell}$ contains as a subgraph any graph on $\ell$ vertices, we get that the statement of Claim~\ref{sequence} is also true for any graph $Y$.\\

Below, we use a notion of $(X,Y)$-extension introduced by Spencer in~\cite{Spencer_extensions}. Let $x\in\mathbb{N}$ and $X=\{\omega_1,\ldots,\omega_x\}$ be a set of $x$ vertices called {\it roots}. Let $Y$ be a graph on $\{\omega_1,\ldots,\omega_y\}$, $y>x$. Then a graph $\tilde Y$ on $\{\tilde\omega_1,\ldots,\tilde\omega_y\}$ is called {\it $(X,Y)$-extension of $\tilde X=\{\tilde\omega_1,\ldots,\tilde\omega_x\}$}, if, for distinct $i\in[y]$, $j\in[y]\setminus[x]$, the presence of the edge $\{\omega_i,\omega_j\}$ in $Y$ implies the presence of the edge $\{\tilde\omega_i,\tilde\omega_j\}$ in $\tilde Y$.

In~\cite{Spencer_extensions}, it is proven (by a straightforward application of another Janson's inequality,~\cite[Theorem 2.18 (i)]{Janson_book}) that $[x]$ does not have an $(X,Y)$-extension with probability at most $e^{-\Omega(n)}$. By the union bound, this observation implies the following.

\begin{claim}
Let $\varepsilon>0$, $\tilde n\in(\varepsilon n,n]$ be a sequence of positive integers. Then whp 
\begin{itemize}
\item for any pair $u,v\notin[\tilde n]$ of adjacent in $G(n,p)$ vertices, there exists a copy $H_{uv}$ of $H'$ in $G(n,p)|_{[\tilde n]}$ such that $(v,u,H_{uv})$ is $H$-completable in $G(n,p)$;
\item for any $v\in[n]\setminus[\tilde n]$, there exists a copy $H_v$ of $H'$ in $G(n,p)|_{[\tilde n]}$ such that $(v,H_v)$ is $H$-completable in $G(n,p)$.
\label{pairs}
\end{itemize}
\end{claim}

Let $b\in\mathbb{N}$. The number of common neighbors of $[b]$ in $G(n,p)$ has binomial distribution with parameters $n-b$ and $p^b$. By the Chernoff bound, this number is smaller than $\frac{1}{2}p^b(n-b)$ with probability at most $e^{-\Omega(n)}$. By the union bound, we get the following.

\begin{claim}
Let $\varepsilon>0$, $\tilde n\in(\varepsilon n,n]$ be a sequence of positive integers, $b\in\mathbb{N}$. Then whp, any subset of $[n]\setminus[\tilde n]$ of size at most $b$ has at least $\frac{\varepsilon}{2}p^b n$ common neighbors in $G(n,p)|_{[\tilde n]}$.

%Let $A, B, A\cap B = \emptyset$ be subsets of $V(G(n, p))$ such that edges between $A$ and $B$ do not depend on the choice of $A$ and $B$ and whp $A = \Omega(n)$.  Let $Y = \{y_1, \ldots y_j\}$ be a set of subsets of $B$, such that the largest of them has size $C_1 = O(1)$, $j = j(n)$.

%Let $A_i$ be the set of common neighbours of $y_i$ in $A$, $1\le i \le j$.

%Then whp $|A_i| = \Omega(n)$ for all $1 \le i \le j$.
\label{simple}
\end{claim}

%{\it Proof}

%Notice that there are at most $O(n^{C_0})$ disticnt subsets of $B$ of size not greater than $C_1$, so we can assume that $j(n) = O(n^{C_0})$.

%Fix some $i, 1\le i\le j$. Then for every $u \in A$ the probability that $u$ is common neighbour of $y_i$ is $p^{|y_i|} \ge p^{C_1}$. So, if $|A| = a$ then $|A_i| \ge \xi_i$ where $\xi_i$ has distribution $\Bin(a, p^{C_1})$.

%By Chernoff's inequality 

%$${\sf P}(|A_i| \le a p^{C_1} / 4||A| = a) \le e^{-\Omega(a)}$$

%Notice that whp $|A| = \Omega(n)$.

%Then 
%$${\sf P}\forall i: |A_i| = \Omega(n)) \ge 1 - o(1) - j(n)e^{-\Omega(n)} \to 1,  n \to \infty.$$

%Then whp $|A_i| = \Omega(n)$ for all $1 \le i \le j$.  \qed

Now, let us prove that there exist $c>0$ and sets $K,S,D,Z_1,Z_2,Z_3,
T,V^1,V^2,V^3$ such that 

\begin{itemize}
\item[---] the tuple $(G(n, p),c,K,S,D,Z_1,Z_2,Z_3,T)$  whp satisfies the first $H$-saturation property, 

\item[---] the tuple $(G(n,p),V^1,V^2,V^3,S,D)$ whp satisfies the second $H$-saturation property.
\end{itemize}

%Notice that if some two properties hold whp then whp they hold at the same time.

Let us start with the first $H$-saturation property. We will at the same time define the parameters and prove that whp each condition (out of 7 from the definition of the first property) holds for these parameters. %Then the first $H$-saturation property holds for these parameters.

\begin{enumerate}
    \item[1.] Since $G_1 \overset{d}{=} G(n/2, p)$, whp, in $G_1$, there is a clique of size at least $c \ln n$ for some positive constant $c$ (see ~\cite{Janson_book}, Theorem 7.1). Let $K$ be this clique of size $k \ge  c \ln n$. So, Condition \ref{K} of the first property holds whp.
    
\end{enumerate}

Let 

\begin{itemize}
\item $S$ be a set of $|S^0_k|$ vertices of $K$, 
\item $D$ be a set of $d$ vertices of $S$,
\item $Z$ be the set of all common neighbours of $D$ from $V(G_2)$.
\end{itemize}

Notice that edges between $D$ and $V(G_2)$ do not depend on the choice of $K, S$ and $D$. Then, $|Z|$ has binomial distribution with parameters $n/2$ and $p^d$. Therefore, whp $\frac{1}{4}p^d n < |Z| < \frac{3}{4} p^d n$ (say, by Chebyshev's inequality).

\begin{enumerate}
    \item[2.] Let $R$ be an arbitrary set of $r$ vertices in $K \setminus S$. %Then edges between $R$ and $Z$ do not depend on the choice of $K, S, D, Z$ and therefore have independent Bernoulli distributions with parameter $p$.
    Then $T$ is the set of all common neighbours of $D \sqcup R$ in $V(G_2)$.

    Edges between $V(G_2)$ and $D \sqcup R$ do not depend on the choice of $D$ and $R$. Then, $|T|$ has binomial distribution with parameters $n/2$ and $p^{r + d}$. Therefore, whp $\frac{1}{4}p^{r+d} n < |T| < \frac{3}{4} p^{r+d} n$.% (say, by Chebyshev's inequality).%by Chernoff's inequality

    %$${\sf P}(|T| \le  p^{r + d} n/4) \le e^{-\Theta(n)}.$$
    
    %and whp  $|T| = \Omega(n)$.
    
    %Fix $v \in Z \setminus T$.
    %Remember that $U_v = T \cap N_G(v)$ for every $v \in Z \setminus T$.
    
    %Notice that edges between $Z \setminus T$ and $T$ do not depend on the choice of $Z$ and $T$.
    
    %Suppose that $|T| = t$. Then $|U_v|$ is distributed as $\mathsf{Bin}(t, p)$.
    
    %So, by Chernoff's bound
    
    %$${\sf P}(|U_v| \le  p t/2 \ |\ |T| = t) \le e^{-\Theta(t)}.$$
    
    %$${\sf P}(\exists v \in Z \setminus T: |U_v| \le  p t/2 \ |\ |Z| = z, |T| = t) \le (z-t)e^{-\Theta(t)} \le n e^{-\Theta(t)}.$$

    %$${\sf P}(\exists v \in Z \setminus T: |U_v| \le  p \gamma_1 n/8) \le o(1) +  ne^{-\Theta(n)} \to 0.$$

    %So, whp $|U_v| = \Omega(n)$.
    
    Notice that edges between $T$ and $Z\setminus T$ and edges inside $T$ are independent of the choice of $Z$ and $T$, so, conditioned on $Z$ and $T$, they have independent Bernoulli distributions. Then, by Claim \ref{pairs},  whp for every adjacent in $G(n, p)$ pair $v_1, v_2 \in Z\setminus T$ there exists a copy $H_{v_1, v_2} \cong H'$ inside $G(n, p)|_T$ such that $(v_1, v_2, H_{v_1, v_2})$ is $H$-completable. %which is the condition \ref{Z_T} of the first property.
    % \item[2)] 

%{\bf ???} It implies that whp $|Z| > \gamma n/2 + \varepsilon/2$.

%Let $c$ be less than $\gamma n/2 + \varepsilon/2$. Then whp $|Z| \ge cn$.
%
%Notice that the edges between  $[n] \setminus (K\sqcup Z)$ and $Z$ are independent of the choice of $K$, $S, D$ and $Z$, so, conditioned on $K$, $S, D$ and $Z$, they are distributed as independent Bernoulli variables. Therefore, by the Chernoff bound and the union bound, whp every $v \in [n] \setminus (K\cup Z)$ has at least $d$ neighbors in $Z$, as the probability of the complement is at most
%$$
% n{\sf P}(\mathrm{Bin}(\gamma n/4,p)<d)= %? ????? ????? ????????
% n\exp(-\Omega(n))\to 0,\quad n\to\infty.
%$$

%So, condiditon \ref{construction_need}  of the first property whp holds.

\item[3.] Let $Z = Z_1 \sqcup Z_2 \sqcup Z_3$ be an equal partition of $Z$ (the sizes of the parts differ by at most one). Then whp $|Z_1| = \Omega(n)$ and edges indside $Z_1$ and between $Z_1$ and $[n] \setminus (K \sqcup Z)$ do not depend on the choice of $K, Z$ and  $Z_1$. Therefore, by Claim \ref{pairs}, for  whp for every $v \notin K \sqcup Z$,  there exists  $H_v \cong H'$ in $G(n, p)|_{Z_1}$ such that $(v, H_v)$ is $H$-completable. It implies that there is some copy $\hat H \cong H|_{V(H)\setminus \{v\}}$ such that $w_1$ is mapped onto $v$ and $H'$ is mapped onto $H_v$. Clearly, it implies the existence of $V_v\subset V(H_v)\setminus\{v\}$ (by the definition) such that $(v, H_v, V_v)$ is $H$-completable. % Let us define $V_v = N_{\hat H}(v)$. Then 
%Then by defining $V_v$ as the set of the second ends of the edges from $v$ in some  copy of $H|_{V(H) \setminus \{w_2\}}$ on $\{v\} \cup H_v$ such that $w_1$ is mapped onto  $v$ and $H'$ onto $H_v$ (such copy exists by definition of $H$-completable tuple)
%we get that $(v, H_v, V_v)$ is $H$-completable which implies condition \ref{H_v} of the first property.

\item [4.] Recall that, for a vertex $v \notin K \sqcup Z$, we denote by $U_v(H_v, V_v)$ the set of all neighbours $u$ of $v$ in $Z_2$ such that the tuple $(v,u,H_v,V_v)$ is $H$-completable.

Set $\tilde U_v: = U_v(H_v, V_v)$. Fix $v_1, v_2 \notin K\sqcup Z$. Notice that if $u\in Z_2$ is a common neighbour of $V(H_{v_1}) \cup V(H_{v_2})\cup\{v_1,v_2\}$ then $u$ lies in $\tilde U_{v_1} \cap \tilde U_{v_2}$. Notice that $V(H_{v_1}) \cup V(H_{v_2}) \subset Z_1$ and $|V(H_{v_1}) \cup V(H_{v_2})| \le 2r$. Edges between vertices of $[n]\setminus(K\sqcup Z)$ and vertices of $Z_2$ are independent on the choice of these sets. Edges between $Z_1$ and $Z_2$ are independent of the choice of $Z_1$ and $Z_2$ as well, and, moreover, $Z_2 = \Omega(n)$. So, applying Claim \ref{simple}, we get that whp $|\tilde U_{v_1} \cap \tilde U_{v_2}| = \Omega(n)$.  %there are $\Omega(n)$ common neighbours in $Z_2$ of $V(H_{v_1}) \cup V(H_{v_2})\cup\{v_1,v_2\}$.  
%Then whp $|\tilde U_{v_1} \cap \tilde U_{v_2}| = \Omega(n)$ for all pairs $v_1, v_2 \notin K \sqcup Z$.
%Notice that edges inside $\tilde U_{v_1} \cap \tilde U_{v_2}$ do not depend on the choice of this set, so they still have independent Bernoulli distributions. So, 
By Claim \ref{sequence}, %(there are at most $n^2$ different sets $\tilde U_{v_1} \cap \tilde U_{v_2}\subset Z_2$, $v_1,v_2\notin K\sqcup Z$), 
whp there is a copy of $H'$ in every $\tilde U_{v_1} \cap \tilde U_{v_2}$.
    
\item [5.] Fix $v \notin K \sqcup Z$, $u \in (K \setminus S)\sqcup Z_2\sqcup Z_3$. Notice that edges between $([n] \setminus K \sqcup Z) \sqcup (K \setminus S)\sqcup Z_2  \sqcup Z_3 = [n] \setminus S \setminus Z_1$ and $Z_2$ are independent of the choice of $S, K, Z, Z_1, Z_2, Z_3$ and so are edges between $Z_1$ and $Z_2$. %So, they are distributed as independent Bernoulli variables. 
By the Chernoff bound, with probability $1-e^{-\Omega(n)}$, vertices from $V(H_{v}) \cup \{u,v\}$ have $\Omega(n)$ common neighbors in $Z_2$. Since there are at most $n^2$ choices of $v \notin K \sqcup Z$, $u \in (K \setminus S)  \sqcup Z_2\sqcup Z_3$, by the union bound, we get that whp there are $\Omega(n)$ common neighbours in $Z_2$ for each element of $\{V(H_{v}) \cup \{u,v\} \ | \ v \notin K \sqcup Z, u \in (K \setminus S)  \sqcup Z_2\sqcup Z_3\}$. Notice that, if $u' \in Z_2$ is common neighbour of $V(H_v) \cup \{u,v\}$, then $u' \in \tilde U_v \cap N_{G(n, p)}(u)$. So, whp, for every $v \notin K \sqcup Z$, $u \in (K \setminus S)\sqcup Z_2 \sqcup Z_3$,   $|\tilde U_v \cap N_{G(n, p)}(u)| = \Omega(n)$. %Edges inside $\tilde U_v \cap N_{G(n, p)}(u)$ do not depend on the choice of $\tilde U_v$ and $u$ and there are at most $n^2$ pairs $(u, v)$. 
By Claim \ref{sequence}, whp there is a copy of $H'$ in every $\tilde U_v \cap N_{G(n, p)}(u)$.

%\textcolor{red}{{ WHAT IF $u \in Z_2$ }}

\item [6.] Notice that edges between $([n] \setminus K \setminus Z) \sqcup Z_1$ and $Z_3$ do not depend on the choice of these sets, and whp $|Z_3|= \Omega(n)$,   therefore whp any $v \notin K \sqcup Z, u \in Z_1$ have $\Omega(n)$ common neighbours in $Z_3$ by Claim \ref{simple}. %Edges inside each $Z_3 \cap N_{G(n, p)}(u) \cap N_{G(n, p)}(v)$ do not depend on the choice of the above sets and 
%There are at most $n^2$ choices of $(u,v)$. 
So, by Claim \ref{sequence}, whp there is a copy of $H'$ in  $Z_3 \cap N_{G(n, p)}(u) \cap N_{G(n, p)}(v)$ for every $u \in Z_1$, $v \notin K \sqcup Z$.

\item [7.] Notice that, if $u \notin S$, then edges between $u$ and $Z$ do not depend on the choice of $Z$ and $u$. If $u \notin S$ then, by the Chernoff bound, with probability $1-e^{-\Omega(n)}$, $u$ has $\Omega(n)$ neighbors in $Z$. Since there are at most $n$ choices of $u \notin S$, by the union bound, we get that whp, every vertex $u \notin S$ has at least $\Omega(n)$ neighbors in $Z$. 
%
%by Claim \ref{simple}, whp there are $\Omega(n)$ neighbours of $u$ in $Z$ (in particular, when $u\in Z$, we do not include $u$ in $[\tilde n]$ in notations from Claim~\ref{simple}). Notice that edges inside $N_G(u) \cap Z$ do not depend on the choice of these sets and there are at most $n$ choices of $U$. 
So, we can apply Claim \ref{sequence} and get that whp there is a copy of $H'$ in  $N_{G(n, p)}(u) \cap Z$ for every $u \notin S$.

%\textcolor{red}{LEFT $u \in Z$}

\end{enumerate}

Now let us prove the second property. Let $V(G_2)=V^1\sqcup V^2\sqcup V^3$ be an equal partition.

\begin{enumerate}
    \item Let us find a copy of $H'$ such that $N_H(w_1) \setminus \{w_2\}$ is mapped onto $D$ and other vertices of this copy lie in $V^1$. Notice that edges between $D$ and $V^1$ and inside $V^1$ have independent Bernoulli distributions. Recall that $D$ induces a clique in $G(n,p)$ and $|V^1|=\Omega(n)$. So, Claim~\ref{simple} implies the existence (whp) of $\Omega(n)$ common neighbors of $D$ in $V^1$. By Claim~\ref{sequence}, there exists an $r$-clique in the set of all common neighbors of $D$. This immediately implies the existence of the desired $H_D$.

%    applying the same method as in \cite{Spencer's extensions} we get that there whp exists such a copy.% with probability at most $e^{-\Omega(n)}$.
    
%    Let $H_D$ be this copy. 
    
    Notice that all common neighbours of $V(H_D)\cup\{v\}$ in $V^2$ lie in $\hat U_v$ for every $v \in S \setminus D$ as $D \subset N_{G(n, p)}(v)$. %therefore $(v, u, H_D, D)$ is $H$-completable. 
    Edges between $D \sqcup V^1 \sqcup (V(G_1)\setminus S) \sqcup V^3$ and $V^2$ do not depend on the choice of these sets. So, by Claim \ref{simple}, whp, for every $u\in V(G_1)\setminus S) \sqcup V^3$, there are $\Omega(n)$ common neighbours of $V(H_D) \sqcup \{v,u\}$ in $V^2$ (and so $|\hat U_v \cap N_{G(n,p)}(u)| = \Omega(n)$). By Claim \ref{sequence}, whp, for any $v \in S \setminus D$, $u \in (V(G_1)\setminus S) \sqcup V^3$, there is a copy of $H'$ in $\hat U_v \cap N_{G(n, p)}(u)$.
    
    \item Edges between $(S \setminus D) \sqcup V^1 \sqcup V^2$ and $V^3$ do not depend on the choice of these sets, so, by Claim \ref{simple}, for any $v \in (S \setminus D)$, $u \in V^1 \sqcup V^2$, whp there are $\Omega(n)$ common neighbours of $\{u, v\}$ in $V^3$. By Claim \ref{sequence}, whp, for any $v \in (S \setminus D)$, $u \in V^1 \sqcup V^2$, there is a copy of $H'$ inside $V^3 \cap N_{G(n, p)}(u) \cap N_{G(n, p)}(v)$.
    
\end{enumerate}

%Let us now define $D$. The construction of $F$ depends only on $S, K$ and $Z$ which we have already defined. 
%Notice that in $F$ there are $d$ edges from every vertice of $Z$ to $S$.  

%Define $N_z = N_F(z) \cap S \subset S$ for $z \in Z$. Then $|N_z| = d$ for every $z \in Z$.

%It implies that there exists $D \subset S$ such that $D$ equals $N_z$ for at least  $\frac{|Z|}{C_{|S|}^d}$ vertices $z \in Z$. 

%We know that whp $|Z| \ge cn$, therefore  if we take $c < \frac{1}{{C_{|S|}^d}}$ then $D$ whp equals $N_z$ for $c^2 n$ vertices $z \in Z$ which implies that whp $D$ has at least $c^2 n$ neigbours in $Z$ and first property \ref{D} whp holds.

\section{Proof of Theorem \ref{threshold_stars}}
 Recall that $t\geq 3$ is assumed.\\

Let us first notice that, if $G=G_1\sqcup\ldots\sqcup G_m$ consists of $m$ %un?
connected components $G_1,\ldots,G_m$, then 
$$
\mathrm{wsat}(G,K_{1,t})=\sum_{i=1}^m\mathrm{wsat}(G_i,K_{1,t}).
$$
Therefore, $\mathrm{wsat}(G,K_{1,t})$ is at least the number of non-empty components in $G$. \\

Let $\frac{1}{n^2}\ll p\leq\frac{\ln n}{2n}$. Consider $w_n>0$ such that 
\begin{itemize}
\item $w_n\to\infty$ as $n\to\infty$, 
\item $w_n\leq\frac{\ln n}{2}$ for $n$ large enough,
\item $\frac{w_n}{n^2}\leq p$ for $n$ large enough. 
\end{itemize}
Then, for $n$ large enough, $\frac{w_n}{n}\leq pn\leq\frac{\ln n}{2}$. Let $X$ be the number of isolated edges in $G(n,p)$. Since 
$$
{\sf E}X={n\choose 2}p(1-p)^{2(n-2)}\sim \frac{1}{2}\exp\left[\ln n+\ln (np)-2pn\right]\geq\frac{w_n}{2}(1+o(1))
$$ 
and 
$$
\mathrm{Var}X={\sf E}X+{n\choose 2}{n-2\choose 2}p^2(1-p)^{4+4(n-4)}-({\sf E}X)^2={\sf E}X+O\left(({\sf E}X)^2\frac{\ln n}{n}\right),
$$
by Chebyshev's inequality, we get that
$$
 {\sf P}\left(X< \frac{w_n}{3}\right)\leq\frac{\mathrm{Var}X}{({\sf E}X-w_n/3)^2}\to 0, \quad n\to\infty.
$$
Therefore, whp 
$$
\mathrm{wsat}(G(n,p),K_{1,t})\geq w_n>{t\choose 2}=\mathrm{wsat}(K_n,K_{1,t}).
$$

Now, let $p>\frac{\ln n}{2n}$. For such $p$, for any $\varepsilon>0$, whp there exists a connected component $G_n$ in $G(n,p)$ of size at least $(1-\varepsilon)n$~\cite[Theorem 5.4]{Janson_book}. Clearly, $\mathrm{wsat}(G(n,p),K_{1,t})\geq\mathrm{wsat}(G_n,K_{1,t})$.\\

Below, for a connected graph $G$, we give a necessary and sufficient condition, in terms of the existence of a subgraph from a certain class, for the stability property $\mathrm{wsat}(G,K_{1,t})={t\choose 2}$.

Let $t\leq y<x\leq n$ be integers and $G$ be a graph on $[n]$. Let $F'\subset F\subset G$, $V(F)=\{v_1,\ldots,v_x\}$, $V(F')=\{v_1,\ldots,v_y\}$. Let us call $F$
a \textit{saturating structure of length $x$ in $G$ with the core $F'$}, if every $v_i$, $y+1 \le i \le x$, sends exactly $t-1$ edges to the previous vertices $v_1,\ldots,v_{i-1}$ in $F$, i.e. $|N_F(v_i) \cap \{v_1, \ldots, v_{i-1}\}| = t-1$. We call $y$ {\it the size} of the core. The vector $\mathbf{v}=(v_1,\ldots,v_x)$ is called a {\it saturating ordering of $F$}. %Finally, let us call the pair $(F,\mathbf{v})$ an {\it ordered saturating structure}.

%Let us prove that the existence of a saturating structure is necessary and sufficient for the stability property.\\

\begin{claim}
Let $G$ be connected. 
\begin{enumerate}
\item If $G$ contains a saturating structure of length $n$ with a core $F'\cong K_t$, then $\mathrm{wsat}(G,K_{1,t})={t\choose 2}$.
\item If $\mathrm{wsat}(G,K_{1,t})={t\choose 2}$ and $G$ has at least $\mu$ vertices with degrees at least $t-1$, then $G$ contains a saturating structure of length $\mu$ with a core of size at most ${t+1\choose 2}$.
\end{enumerate}
\label{equiv_condition_sat_structure}
\end{claim}

{\it Proof.} First, assume that $G$ contains a saturating structure $F$ of length $n$ with a core $F' \cong K_t$. It is clear that $F$ is both weakly $(G,K_{1,t})$-saturated and weakly $(K_n,K_{1,t})$-saturated. In particular, it implies that $G$ is weakly $(K_n,K_{1,t})$-saturated and, therefore, $\mathrm{wsat}(G,K_{1,t})\geq\mathrm{wsat}(K_n,K_{1,t})$. Since $F$ is weakly $(G,K_{1,t})$-saturated, it remains to prove that $\mathrm{wsat}(F,K_{1,t})\leq{t\choose 2}$. Let $(v_1,\ldots,v_n)$ be a saturating ordering of $F$. Then $F'=F|_{\{v_1,\ldots,v_t\}}\cong K_t$ has exactly $t\choose 2$ edges. Let us show that $F'$ is weakly $(F,K_{1,t})$-saturated. Since each vertex of $v_1,\ldots,v_t$ has degree $t-1$ in $F'$, we can restore all edges in $F$ adjacent to one of these vertex. In particular, we restore all edges going from $v_{t+1}$ to $v_1,\ldots,v_t$. Proceeding in this way by induction, we restore all edges of $F$.\\

%    If $F' = F|{\{v_1, \ldots v_t\}}$ then $F' \cong K_t$ and we can restore all edges of $F$ going out of $\{v_1, \ldots v_t\}$. Then, there are $t-1$ edges from $v_{t+1}$ restored, so we can restore all other edges going out of $v_{t+1}$. Similarly, on each step we can add $v_i$ and restore all edges of $F$ from $v_i$ for $t + 1 \le i \le n$. Therefore, $\mathrm{wsat}(F, K_{1, t}) \le {t \choose 2}$. 
    
 %   Notice that $F$ is $(K_n, K_{1, t})$-saturated, so $\mathrm{wsat}(F, K_{1, t}) \ge \mathrm{wsat}(K_n, K_{1, t}) = {t \choose 2}$.

Now, let $\mathrm{wsat}(G,K_{1,t})={t\choose 2}$ and $G$ have at least $\mu$ vertices with degrees at least $t-1$. Let $F'$ be a weakly  $(G,K_{1,t})$-saturated graph with ${t\choose 2}$ edges and $y$ non-isolated vertices. Let us order these vertices of $F'$ in a way $v_1,\ldots,v_y$ such that $v_i$ plays the role of the $i$th central vertex of $K_{1,t}$ in a $K_{1,t}$-bootstrap percolation process that starts on $F'$ and finishes on $G$. Clearly, for every $i\in[t-1]$, the vertex $v_i$ sends at least $t-i$ edges to $F'|_{\{v_{i+1},\ldots,v_y\}}$. Since the total number of these edges is ${t\choose 2}$, $F'$ can not contain any other edge. The bound $y\leq t+{t\choose 2}={t+1\choose 2}$ follows.

Consider a $K_{1,t}$-bootstrap pecolation process that starts on $F'$ and finishes on $G$. Let $e_1,\ldots,e_m$ be edges appearing in this process sequentially that contain at least one vertex outside $\{v_1,\ldots,v_y\}$. Let $w_1,\ldots,w_{\mu-y}$ be vertices of $G$ outside $\{v_1,\ldots,v_y\}$ with degrees at least $t-1$ ordered in the following way. 

\begin{itemize}
\item Let $i\in[m]$ be such that $e_i$ contains $w_1$ and a vertex from $\{v_1,\ldots,v_y\}$, there are exactly $t-2$ edges among $e_1,\ldots,e_{i-1}$ that contain $w_1$, and all of them have the second end in $\{v_1,\ldots,v_y\}$. 
\item For $j\in\{2,\ldots,\mu-y\}$, let $i_j\in[m]$ be such that $e_{i_j}$ contains $w_j$ and a vertex from $\{v_1,\ldots,v_y,w_1,\ldots,w_{j-1}\}$, there are exactly $t-2$ edges among $e_1,\ldots,e_{i_j-1}$ that contain $w_j$, and all of them have the second end in $\{v_1,\ldots,v_y,w_1,\ldots,w_{j-1}\}$. 
\end{itemize}

Such an ordering exists due to the definition of the  $K_{1,t}$-bootstrap pecolation process. Then, the desired saturating structure of length $\mu$ is obtained from $F'|_{\{v_1,\ldots,v_y\}}$ by adding $w_i$, $i\in[\mu-y]$, with the $t-1$ edges going to the previous vertices $v_1,\ldots,v_y,w_1,\ldots,w_{i-1}.$ $\Box$\\

Now, due to Claim~\ref{equiv_condition_sat_structure} and the fact that whp $G(n,p)$ has at least $n/2$ vertices with degrees at least $t-1$ (this can be proven by a straigtforward application of Chebyshev's inequality to the number of vertices with such degrees), Theorem~\ref{threshold_stars} immedialtely follows from

\begin{claim}
\begin{enumerate}
    \item There exists $c>0$ such that, if $p<cp(n, t)$, then whp there is no saturating structure of length $\lfloor\ln n\rfloor$ and with a core of size at most ${t+1\choose 2}$ in $G(n, p)$.
    \item There exists $C>0$ such that, if $p> Cp(n, t)$, then  %If $p\gg n^{-\frac{1}{t-1}}$, then 
    whp there exists a saturating structure of length $n$ with a core isomorphic to $K_t$ in $G(n, p)$. 
\end{enumerate}

\label{existence_saturating}
\end{claim}

\subsection{Proof of Claim~\ref{existence_saturating}.1}

%Denote by $\gamma_n = \ln n^{\frac{1}{s-1}}

Let $x=\lfloor\ln n\rfloor$, $y={t+1\choose 2}$, $c<e^{-(y+1)/(t-1)}$. Let $p<cp(n,t)$.

Let $X$ be the number of subgraphs $F$ in $G(n,p)$ on $x$ vertices such that there exist $v_1,\ldots,v_i\in V(F)$, $i\in\{t-1,t,\ldots,y\}$ satisfying the following property:
\begin{center}
 for every $v\in V(F)\setminus\{v_1,\ldots,v_i\}$, there exists a set $N_v$ of its $t-1$ neighbors in $F$ such that, for every $u\in N_v$, $v\notin N_u$.
\end{center}

Clearly, an ordered saturating structure of size $x$ with a core of size at most $y$ is a subgraph with the above property. Therefore, it is sufficient to prove that ${\sf P}(X\geq 1)\to 0$ as $n\to\infty$.

Let us bound from above ${\sf E}X$:
$$
 {\sf E}X\leq \sum_{i=t-1}^{y}{n\choose x}{x\choose i}{x\choose t-1}^{x-i}p^{(x-i)(t-1)}\leq \frac{n^xx^y}{x!}\sum_{i=t-1}^y \left({x\choose t-1}p^{t-1}\right)^{x-i}
$$
$$
 \leq \frac{yn^xx^y}{x!}\left(xp\right)^{(t-1)(x-y)}=e^{x\ln n-x\ln x+x+(x-y)(t-1)\ln (xp)+O(\ln\ln n)}.
$$
Since 
$$
\ln p<\ln (cp(n,t))=-\frac{1}{t-1}\ln n-\frac{t-2}{t-1}\ln\ln n-\ln\frac{1}{c},
$$
we get that
$$
 {\sf E}X\leq\exp\left[y\ln n-y\ln\ln n+x-(x-y)(t-1)\ln\frac{1}{c}+O(\ln\ln n)\right]=
$$
$$
 \exp\left[\ln n\left(y+1-(t-1)\ln\frac{1}{c}\right)+O(\ln\ln n)\right]\to 0,\quad n\to\infty.
$$
Markov's inequality implies ${\sf P}(X\geq 1)\to 0$.

\subsection{Proof of Claim~\ref{existence_saturating}.2}

%Let us recall that the {\it $k$th power of a Hamiltonian path} is a graph with the vertex set $\{v_1,\ldots,v_n\}$ and edges $\{v_i,v_{i+j}\}$ for all $j\in[k]$ and all $i\in[n-k]$. Clearly, if a graph $G$ on the vertex set $[n]$ contains the $(t-1)$th power of a Hamiltonian path, then it also contains a saturating structure of length $n$ with a core isommorphic to $K_t$. Therefore, for such $G$, $\mathrm{wsat}(G,K_{1,t})={t\choose 2}$. Then, Claim~\ref{existence_saturating}.2 immediately follows from 

%\begin{theorem}[K\"{u}hn, Osthus, 2012, \cite{KO}; Kahn, Narayanan, Park, 2020, \cite{Kahn}]
%Let $k\geq 2$. If $p\gg n^{-1/k}$, then whp $G(n,p)$ contains the $k$th power of a Hamiltonian path. 
%\end{theorem}

Set $p=Cp(n, t)$ where $C$ is a large positive constant (for example, any $C$ bigger than $2^{\frac{2t-1}{t-1}}t(t-1)$ is sufficient). Since `containing a saturating structure of length $n$ with a core isomorphic to $K_t$' is an increasing property, it is sufficient to prove that it holds whp for this value of $p$.

%We will define $C$ later.

%Set $x = \lfloor \ln n \rfloor$.

The structure of the proof is the following: at first we prove the existence of a saturating structure of size $x = \lfloor \ln n \rfloor$ whp, then we extend this structure to size $y = \left\lfloor{\left(\ln n\right)^{\frac{t-1}{t-2}}}\right\rfloor$  and, finally, we extended it to the desired size $n$.

\subsubsection{Saturating structure of size $x = \lfloor \ln n \rfloor$}

Let $X$ be the number of saturating structures $F$ of size $x$ in $G(n, p)$ with a core isomorphic to $K_t$ and a saturating ordering $(v_1, \ldots, v_x)$ such that $v_1  < v_2 < \ldots < v_x$. Let us call such an ordering a {\it canonical} saturating ordering. Let $S$ be the set of all such structures in $K_n$. Then $X = \sum\limits_{A \in S} I_A$ where $I_A$ indicates that $A$ belongs to $G(n, p)$. We have 
\begin{equation}
\Exp X \sim \frac{n^x}{x!} p^{{t-1 \choose 2}} \prod_{i=t}^{x}\left[{i-1 \choose t-1} p^{t-1}\right]. 
\label{expectation_asymp}
\end{equation}

Notice that  
$$
\prod\limits_{i=t}^{x}{i-1 \choose t-1} = \frac{1}{[(t-1)!]^{x-t}} \frac{t!\cdot\ldots\cdot(x-1)!}{1!\cdot\ldots\cdot(x-t)!} =
$$ 
$$ 
\frac{1}{[(t-1)!]^{x-t}} \frac{(x-t+1)!\cdot\ldots\cdot(x-1)!}{1!\cdot\ldots\cdot(t-1)!} \ge  \frac{[(x-t+1)!]^{t-1}}{[(t-1)!]^{x-t}} \sim $$ $$ \sim \frac{(2\pi x)^{(t-1)/2}  (x-t+1)^{(x-t+1)(t-1)}}{e^{(x-t+1)(t-1)}[(t-1)!]^{x-t}}>
(x-t+1)^{(x-t+1)(t-1)}.
$$

%For large $n$ this expression is not less than $(x-t+1)^{(x-t+1)(t-1)}.$
%\\ \sim \frac{n^x}{\tqrt{2\pi x}(x/e)^x} $  $p^{{t-1 \choose 2} + (t-1)(x-t+1)} \prod_{i=t}^{x}{i-1 \choose t-1} 

So, 
$$
\Exp X  \ge (1+o(1)) \frac{n^x}{\sqrt{2\pi x}(x/e)^x}  (x-t+1)^{(t-1)(x-t+1)}  p^{{t-1 \choose 2} + {(t-1)(x-t+1)}} = 
$$
$$
\exp \left\{ %x \ln n - \ln \sqrt{2\pi x} - x \ln x + x - x\ln A + (t-1)(x-t+1)\ln (x-t+1) + (x - t/2)(t-1) \ln p =
\frac{t}{2} \ln n + x \left[1 + (t-2)\ln x + (t-1) \ln C - (t-2) \ln \ln n\right] + o(x) \right \} = 
$$

%Set $x_0 = B\lfloor \ln n \rfloor$, $B \in \mathbb{R}$ (we will select $B$ later).

%Then $$\Exp X(x_0) 
$$
\exp \left \{ %x \ln n - \ln \sqrt{2\pi x} - x \ln x + x - x\ln A + (t-1)(x-t+1)\ln (x-t+1) + (x - t/2)(t-1) \ln p =
\ln n \left[\frac{t}{2} +  (1 + (t-1) \ln C) \right] + o(\ln n) \right \}.
$$

As $C > e^{-\frac{1}{t-1}}$, we get $\Exp X \to \infty$.\\

Now let us estimate $\Var X$. Since, for disjoint $A,B\in S$, $I_A,I_B$ are independent, we get
$$
\mathrm{Var} X  = \Exp X^2 - (\Exp X)^2 =
$$
\begin{equation}
\Exp \left[\sum\limits_{\substack{A, B \in S,\\ V(A) \cap V(B) = \varnothing}}I_A I_B\right] + 
\Exp\left[ \sum\limits_{\substack{A, B \in S,\\ V(A) \cap V(B) \neq \varnothing}}I_A I_B\right] - 
(\Exp X)^2\leq 
\sum\limits_{\substack{A, B \in S,\\ V(A) \cap V(B) \neq \varnothing}}{\sf E}I_A I_B.
\label{variance:small_sat_structure}
\end{equation}

%\begin{itemize}%[label=$\bullet$]
   % \item 
   
 %  $\Exp \sum\limits_{\substack{A, B \in S,\\ A \cap B = \emptyset}}I_A I_B = (1 + o(1)) \left(\frac{n^x}{x!} \prod_{i=t}^{x}{i-1 \choose t-1} p^{{t-1 \choose 2} + (t-1)(x-t+1)}\right)^2 % \frac{n(n-1)(n-2)...(n-x+1)}{x!} \prod_{i=s}^{x}{i-1 \choose s-1}$ $p^{{s-1 \choose 2} + (s-1)(x-s+1)}$  $\frac{(n-x)(n-x-1)(n-x-2)...(n-2x+1)}{x!}$ $\prod_{i=s}^{x}{i-1 \choose s-1}$ $p^{{s-1 \choose 2} + (s-1)(x-s+1)}$=$
 %   = \\
 %   = (1 + o(1))(\Exp X)^2 $

 %   \item 
    
    Let $A,B\in S$ have a non-empty intersection $W=V(A)\cap V(B)$ and let $w_1<\ldots<w_d$ be the vertices of $W$. %Denote
    
%    Now consider some intersecting structures $A$ and $B$.  Suppose that they have $d$ common vertices $W = \{w_1 < \ldots < w_d\}$. %Notice that canonical saturating orderings of $A$ and $B$ are ordered, so $w_1, \ldots w_d$ li???????, ??? ? ???? ??????????????? ?????? ? ????? ???????????????????, ??????????? ??????? ? B ???? ? ??? ?? ???????, ??? ? ? A. 

%    Common edges can only touch vertices of $W$. Notice that the maximum number of common edges appears when all edges from $W$ go to $W$. 

Since the number of edges in $W$ is maximum when each $w_i$ sends all $\min\{t-1,i-1\}$ edges to the previous $w_1,\ldots,w_{i-1}$, we get that $A$ and $B$ have at most 
$$
M: =\left[{t-1 \choose 2}+ (d-t+1)(t-1)\right]I(d\geq t)+ {d \choose 2}I(d<t).
$$ 
common edges. Notice that $M=(d - t/2)(t-1)$ when $d\geq t$. Denote by $\mathrm{cnt}(d, m_1, m_2)$ the number of pairs $A, B \in S$ such that $|V(A)\cap V(B)|=d$, $A\cap B$ has $m_1$ edges inside the core of $B$ and $m_2$ edges outside the core of $B$. Then 
    \begin{equation}
    \sum\limits_{\substack{A, B \in S,\\ A \cap B \neq \emptyset}}\Exp I_A I_B  =\sum\limits_d \sum\limits_{m_1} \sum\limits_{m_2}  \mathrm{cnt}(d, m_1, m_2) p^{2  z - m_1 - m_2},
    \label{cnt_def}
    \end{equation}
    where
    $$
    z =  {t-1 \choose 2} + (x-t+1)(t-1)  = (x - t/2)(t-1).
    $$

    Notice that 
    \begin{equation}
    \mathrm{cnt}(d, m_1, m_2) \le 
    \frac{EX}{p^z}  {x \choose d} \frac{n^{x-d}}{(x-d)!} 
    \max\limits_{j_{t+1} , \ldots, j_x \in J(d, m_2)}
    \prod_{i=t+1}^{x}{i-1 \choose t-1-j_i} {t-1 \choose j_i},
    \label{cnt_above}
    \end{equation}
    where $J(d, m_2)$ is the set of all tuples $(j_{t+1},\ldots,j_x)\in\{0,1,\ldots,t-1\}^{t-x}$ such that $j_{t+1} + \ldots+ j_x  = m_2$ and the number of non-zero $j_i$ is at most $d$. Clearly, for a $(j_{t+1},\ldots,j_x)\in J(d,m_2)$, we have $\prod\limits_{i=t+1}^{x}{t-1 \choose j_i} \le 2 ^{(t-1)d}$. Moreover,
    \begin{align*}     
    \prod_{i=t+1}^x {i-1 \choose t-1-j_i} & =
    \prod_{i=t+1}^x \frac{(t-1)!(i-t)!}{(t-1-j_i)!(i-t+j_i)!}{i-1\choose t-1} \\ 
     & \leq\prod_{i=t+1}^x\left(\frac{t-1}{i-t+1}\right)^{j_i}{i-1\choose t-1}=
     (t-1)^{m_2}\prod_{i=t+1}^x\frac{{i-1\choose t-1}}{(i-t+1)^{j_i}}.
    \end{align*}

    %\\ $ = O\left(\frac{EX}{p^z} {x \choose d} \frac{n^{x-d}}{(x-d)!} \max_{j_1 +...+ j_x = m_2}\prod_{i=t}^{x}{i-1 \choose t-1-j_i}\right)$
    
  %  Let $k_i = t-1-j_i$. Then
    
    %\frac{cnt(d, m_1, m_2)}{(EX)^2} = 

%    $$\frac{{i \choose k_i}}{{i \choose t-1}} = \frac{(t-1)!(i-t+1)!}{k_i! (i-k_i)!} = \frac{(k_i+1)\ldots(t-1)}{(i-k_i)\ldots(i-t+2)} \le  \left (\frac{t-1}{(i-t+1)}\right )^{j_i}.$$ %where $D = (t-1)^{% = O\left(\frac{1}{i}^{j_i}\right)$.
    
%    So, $\prod\limits_{i=t}^{x}{i-1 \choose t-1-j_i} \le (t-1)^{m_2} \prod\limits_{i=t}^{x} \frac{1}{i-t+1}^{j_i} \prod_{i=t}^{x} {i-1 \choose t-1}$.% for some constant $D>1$.
    
   The function 
   $g(j_{t+1},\ldots,j_x)=\prod\limits_{i=t+1}^{x} \frac{1}{(i-t+1)}^{j_i}$ defined on the intersection of $\{0,1,\ldots,t-1\}^{x-t}$ with the hyperplane $j_{t+1} + \ldots+ j_x  = m_2$ achieves its maximum when 
    $$
    j_i=\left\{
    \begin{array}{cc} 
    t-1, & t+1\leq i\leq t+\left\lfloor \frac{m_2}{t-1}\right\rfloor, \\
    m_2\mod t-1, & i= t+\left\lfloor \frac{m_2}{t-1}\right\rfloor +1, \\ 
    0, & i>t+\left\lfloor\frac{m_2}{t-1}\right\rfloor +1
    \end{array}\right.
    $$
    since $-\ln g=\sum_{i=t+1}^x\alpha_ij_i$ is linear and the coefficients $\alpha_i$ increase as $i$ grows. Therefore, $\prod_{i=t+1}^{x} \frac{1}{(i-t+1)^{j_i}} = O\left[\left(\frac{1}{{\left\lceil\frac{m_2}{t-1}\right\rceil}!}\right)^{t-1} \right]$. Combining this with (\ref{expectation_asymp}) and (\ref{cnt_above}), we get
    \begin{align*}
    f(d, m_1, m_2) & := \frac{\mathrm{cnt}(d, m_1, m_2)p^{2z - m_1 - m_2}}{(EX)^2}  =   \\
    & = O\left({x \choose d} x^d n^{-d}    2 ^{(t-1)d} \frac{(t-1)^{m_2}}{\left(1/\left\lceil\frac{m_2}{t-1}\right\rceil!\right)^{t-1}}  p^{ - m_1 - m_2}\right) =\\ 
    & =O\left(\left(\frac{xe 2 ^{(t-1)} }{d}\right)^d  x^d n^{-d} p^{- m_1}  \left(\frac{(t-1)^2e}{m_2 p}\right)^{m_2} \right).
    \end{align*}
    
 %   Then  $cnt(d, m_1, m_2) = O\left(\frac{EX}{p^z}  {x \choose d} \frac{n^{x-d}}{(x-d)!} (t-1)^{m_2} \left(\frac{1}{{\lceil\frac{m_2}{t-1}\rceil}!}\right)^{t-1} \prod_{i=t}^{x} {i-1 \choose t-1}) \right).$% = O(exp \{d [\ln x + 1 - \ln d] + (x-d) [\ln n + 1 - \ln (x-d)] + m_2[-\ln m_2 -1 + \ln (t-1)]$\ 

    %Set $f(d, m_1, m_2) = \frac{cnt(d, m_1, m_2)p^{2z - m_1 - m_2}}{(EX)^2}$.

%    $$f(d, m_1, m_2) =  %O\left(\frac{{x \choose d} x! n^{x-d} p^{ - m_1 - m_2} }{n^x (x-d)! \prod_{i=t}^{x}{i-1 \choose t-1} }\right) = 
%    O\left({x \choose d} x^d n^{-d}   (t-1)^{m_2}\left(\frac{1}{{\lceil\frac{m_2}{t-1}\rceil}!}\right)^{t-1}  p^{ - m_1 - m_2}\right) = $$ $$\\= O\left(\left(\frac{xe}{d}\right)^d  x^d n^{-d} p^{- m_1} \left(\frac{(t-1)^2e}{m_2 p}\right)^{m_2} \right).$$
    
Since $\left(\frac{(t-1)^2e}{m_2 p}\right)^{m_2}$ increases (as a function of $m_2$) on $(0,M]$ (the maximum is achieved at $\frac{(t-1)^2}{p} \gg M$), we get that 
\begin{equation}
p^{-m_1}\left(\frac{(t-1)^2e}{m_2 p}\right)^{m_2}\leq p^{-m_1} \left(\frac{(t-1)^2e}{(M - m_1) p}\right)^{M - m_1} = p^{-M}\left(\frac{(t-1)^2e}{M-m_1}\right)^{M-m_1}.
\label{last_bound:three_cases}
\end{equation} 
%    Now let us investigate $p^{-m_1} \left(\frac{(t-1)^2e}{(M - m_1) p}\right)^{M - m_1} = p^{-M} \left(\frac{(t-1)^2e}{M - m_1}\right)^{M - m_1} $.

This expression achieves its maximum when $m_1 = M-  (t-1)^2$. Since $m_1\leq{t-1\choose 2}$ by the definition and $M$ may be either large or small depending on the value of $d$, below, we distinguish several scenarios: $d\geq 2t-2$, $t\leq d<2t-2$ and $d<t$.\\

\begin{enumerate}
    \item If $d \ge 2t - 2$, then $M\geq (3t/2-2)(t-1)=(t-1)^2+{t-1\choose 2}$. Therefore, $M-(t-1)^2\geq{t-1\choose 2}$. It means that the bound to the right in (\ref{last_bound:three_cases}) increases with $m_1$ and its maximum value is achieved at $m_1={t-1\choose 2}$. Then,
    \begin{align*}
    f(d, m_1, m_2)   & = O\left(\left(\frac{x^2e 2^ {t-1}}{dn}\right)^d   \left(\frac{(t-1)e}{(d-t+1) }\right)^{(t-1)(d-t + 1)} p^{-M} \right)  = \\  
    & = O\left(e^{d [2\ln x   + 1 - \ln d  - \ln n + (t-1)(\ln 2 - \ln (d - t + 1) +\ln (t-1) +  1 - \ln p)]}\right.\times\\   
    &\quad\times \left. e^{(t/2) (t-1)  \ln p  +  o(x)}\right).
    \end{align*}
    Notice that $\ln (d-t+1) \ge \ln \frac{d}{t}$ (since $d \ge t$). So, 
    \begin{align*}    
    f(d, m_1, m_2)  = & O\left(e^{d\gamma(d)+(t/2) (t-1)  \ln p  +  o(x)}\right),
    \end{align*}
    where
    $$
    \gamma(d)=2\ln x + 1  -  t\ln d  + (t-2) \ln \ln n + (t-1)(\ln t -  \ln C + \ln (t-1) + 1 + \ln 2).
    $$
    Notice that $[d\gamma(d)]' = \gamma(d) - t $ and $\gamma$ decreases. Therefore, $d_0=\gamma^{-1}(t)$ is a point of global maximum of $d\gamma(d)$. Clearly,
    
%    So, if $d_0$ maximizes $dg(d)$ then $g(d_0) = t$. It can be easily seen that such $d_0$ is a point of maximum and not minumum.
    
 %   $$2\ln x  + (t-2) \ln \ln n + (t-1)(\ln t -  \ln C + \ln (t-1) + \ln 2)  = t \ln d_0 .$$
    
    $$
    d_0 = \left(\frac{2(t-1)t}{C}\right)^\frac{t-1}{t}\left(x^2[\ln n]^{t-2}\right)^{\frac{1}{t}} =  \left(\frac{2(t-1)t}{C}\right)^\frac{t-1}{t}x(1+o(1)).
    $$    
    Therefore, $d_0\gamma(d_0)= \left(\frac{2(t-1)t}{C}\right)^\frac{t-1}{t}xt(1+o(1))$.
    As $C > 2^{\frac{2t-1}{t-1}} t(t-1)$, we get that, for $\varepsilon>0$ small enough, $f(d, m_1, m_2)\leq n^{ -\varepsilon}$.
    
\item Let $t \le d < 2t - 2$. Then $\left(\frac{(t-1)^2e}{(M-m_1)}\right)^{M-m_1}\leq e^{(t-1)^2}$. 
Therefore, 
     \begin{align*}
     f(d, m_1, m_2)   & =
     O\left(\left(\frac{x^2}{n}\right)^dp^{-M}\right) \,\,\,\, = 
     O\left(\left(\frac{x^2}{np^{t-1}}\right)^d  
     p^{\frac{t(t-1)}{2}} \right) \\
     & = O\left(\frac{[\ln n]^{td}}{(n [\ln n]^{t-2})^{t/2}}\right)  =  O\left(\frac{1}{n}\right).
     \end{align*}
     % $$  %= O\left(\left(\frac{xe}{d}\right)^d  x^d n^{-d}  \left(\frac{e}{(d - t/2) p}\right)^{(t-1)(d-t/2)} D^x\right) = 
  %  = O(\exp\{\mathbf{d} [2\ln x   + 1  + (t-1) \ln 2 - \ln d  - \ln n   - (t-1) \ln p ] +   t/2 (t-1)  \ln p  +  o(x)\}).$$

%\item  Let $t \le d \le \frac{3t}{2}$.
    
 %   Then $M = {d \choose 2}$.
    
  %  In this case there is maximum when $m_1 = 0, m_2 = (d  - t/2) (t-1)$.
    
   % Again, 
    
    %$$f(d, m_1, m_2)  = %O\left(\left(\frac{x^2e 2^ {t-1}}{dnp^{t-1}}\right)^d  % e^{(t-1)^2} 
     %O\left(\left(\frac{x^2}{np^{t-1}}\right)^d  % e^{(t-1)^2}  p^{t/2(t-1)} \right)  = $$
     
%     $$= O\left(\ln n ^{td} \left(\frac{1}{n \ln n ^{t-2}}\right)^{t/2} \right) = O(n^{-1/3}).$$

\item Finally, let us switch to the case $ d < t$. Since $M-(t-1)^2<0$, we get that $\left(\frac{(t-1)^2e}{M-m_1}\right)^{M-m_1}$ is maximal when $m_1=0$. Therefore,

\begin{align*}
    f(d, m_1, m_2)  & = O\left(\left(\frac{x^2}{n}\right)^d  p^{-M} \right) =
    O\left( \left(\frac{x^2}{np^{\frac{d-1}{2}}}\right)^d\right) \\
    &= O\left( \left(\frac{x^2}{n^{1 - \frac{d-1}{2(t-1)}} [\ln n] ^{\frac{(t-2)(d-1)}{2(t-1)}}}\right)^d\right) = O\left(\frac{\ln^2 n}{n}\right).
\end{align*}
    
%    Then there is maximum when $m_1 = 0, m_2 = {d \choose 2}$.

 %   $f(d, m_1, m_2)  = O\left(\left(\frac{x^2}{n}\right)^d  p^{-{d \choose 2}} \right) = O\left( \left(\frac{x^2}{np^{\frac{d-1}{2}}}\right)^d\right) = \\ = O\left( \left(\frac{x^2}{n^{1 - \frac{d-1}{2(t-1)}} \ln n ^{-\frac{(t-2)(d-1)}{2(t-1)}}}\right)^d\right) = O(n^{-1/3})$
    %O\left(\left(\frac{x^2e}{dn}\right)^d  p^{-{d \choose 2}} \right) = O\left( \left(\frac{x^2e}{dnp^{\frac{d-1}{2}}}\right)^d\right) = \\ = O\left( \left(\frac{x^2}{n^{1 - \frac{d-1}{2(t-1)}} \ln n ^{-\frac{(t-2)(d-1)}{2(t-1)}}}\right)^d\right) \le O(n^{-1/3})$

\end{enumerate}

    %If $x = x_0$ then $d_0 = Dx$.
    %$2 \ln x - t \ln d - (t-2) \ln \ln n + 1 + (t-1) (1 + \ln 2 - \ln C) = t $
    
    %$\frac{2 \ln x - t - (t-2) \ln \ln n + 1 + (t-1) (1 + \ln 2- \ln C)}{t } = \ln d$
    
    %$d_0 \le  \exp \{\frac{2 \ln x - t - (t-2) \ln \ln n + 1 + (t-1) (1 + \ln 2 - \ln C)}{t} \} = O(x^{\frac{-(t-4)}{t}}) << x$.
    
  %  Then $d_0 g(d_0) = D x t$.
    %x = \ln n !!! ???? x = B \ln n
    
%    Then $f(d, m_1, m_2)  =  O(exp\{D x t  +  t/2 (t-1)  \ln p +  o(x)\})$
    
%   If $x = x_0$ then $f(d, m_1, m_2) = O(\exp\{ (D - 1/2)  t \ln n + o (\ln n)\}$.

    %?? B ????? ??????????, ??? ??? ?????? ? ???? ???.

    %Notice that $z = {t-1 \choose 2} + (d-t+1)(t-1)$.
    
    %$\mathbf{d-t/2} [\ln x +1 - \ln d + \ln x - \ln n + (t-1) (1 - \ln (d - t/2)  - \ln p)] \le \mathbf{d-t/2} [\ln x +1 - \ln (d - t/2) + \ln x - \ln n + (t-1) (1 - \ln (d - t/2)  - \ln p)]$

    %Let $f(d, m_1, m_2) = (\frac{xe}{d})^d   (\frac{ne}{(x-d)})^{x-d} D^x (\frac{1}{{\lceil\frac{m_2}{t-1}\rceil}!})^{t-1}$

    %Let  $M(m_2) = \max_{j_1 +...+ j_x = m_2}\prod_{i=t}^{x}C_{i-1}^{t-1-j_i}$

%\end{itemize}
%\vspace{0.2cm}

Combining the above bounds with (\ref{variance:small_sat_structure}) and (\ref{cnt_def}), we get, by Chebyshev's inequality, that, for $n$ large enough,
$$
\mathsf{P} (X = 0) \leq \frac{\Var X}{(\Exp X)^2}  \leq \sum\limits_d \sum\limits_{m_1} \sum\limits_{m_2}  f(d, m_1, m_2) \leq x z^2 n^{-\varepsilon} = o(1).
$$

%Then, by Chebyshev's inequality,

%$$\mathsf{P} (X = 0) \le \mathsf{P}(|\Exp X - X| \ge \frac{\Exp X}{2}) \le \frac{\Var X}{4(\Exp X)^2} \to 0$$

Therefore,  whp there exists a saturating structure of size $x$.%second moment method we get that whp there exists a saturating structure of size $x$. %B???

\subsubsection{Saturating structure of size $y = \left\lfloor{(\ln n)^{\frac{t-1}{t-2}}}\right\rfloor$ }

Let us divide the random graph into two parts: 

$$
G_1 = G(n,p)|_{[\lfloor n/2\rfloor ]},\quad G_2 = G(n, p)|_{[n] \setminus [\lfloor n/2\rfloor ]}.
$$
Let $F_0$ be a saturating structure of size $x = \lfloor\ln\left(\lfloor n/2\rfloor\right)\rfloor$ in $G_1$ (if exists). Let $A_0$ be the event that $F_0$ exists.\\

Let us enlarge the saturating structure by induction.

Let $U_1 \subset V(G_2)$ be a set of vertices connected to at least $t-1$ vertices of $F_0$, for $i=1,2,\ldots$ let $U_{i+1} \subset V(G_2)\setminus (U_1\sqcup\ldots\sqcup U_i)$ be a set of vertices connected to at least $t-1$ vertices of $U_i$. If, for some $\ell$, we get $|V(F_0) \sqcup U_1\sqcup \ldots\sqcup U_{\ell}| \ge y$, then we immediately get a saturating structure of size at least $y$. Let us prove that $\ell = \left\lceil\frac{2}{\ln(t-1)} \ln \ln \ln n\right\rceil$ is the desired value.

Set $X_i = |U_i|$ for $i\in[\ell]$. %Then $|F_0 \cup \ldots F_{\ell}| = \sum\limits_{i=0}^{\ell} X_i$.

Notice that edges between $G_2$ and $F_0$ do not depend on the choice of $F_0$ and have independent Bernoulli distributions. Therefore,
$$
\Exp (X_1|A_0) = \left\lceil\frac{n}{2}\right\rceil P,\quad
\Var (X_1|A_0) = \left\lceil\frac{n}{2}\right\rceil P (1-P),
$$
where
$$
P = \sum\limits_{k = t-1}^{x}{x \choose k} p^k(1-p)^{x-k}.
$$
We get
$$
\Exp (X_1|A_0) \sim \frac{n}{2} {x \choose {t-1}} p^{t-1} (1-p)^{x-t+1} \sim \frac{n}{2} \frac{x^{t-1}}{(t-1)!}p^{t-1} = \frac{C^{t-1}}{2 (t-1)!} x  
$$ 
As $C > (2 (t-1)!)^{\frac{1}{t-1}}$,  we get that $\hat C:= \frac{C^{t-1}}{2 (t-1)!} > 1$. Also, 
$$
\frac{\Var (X_1|A_0)}{(\Exp (X_1|A_0))^2}\sim\frac{1}{{\sf E}(X_1|A_0)} \sim \frac{1}{\hat C x}.
$$
By Chebyshev's inequality, for every $ \delta > 0$, 
\begin{equation}
\mathsf{P}(|X_1 - \hat C x| > \delta \hat C x|A_0) \le (1 + o(1)) \frac{1}{\delta^2 \hat C x}.
\label{A_1st_step}
\end{equation}

Take $\delta>0$ such that $(1-\delta) \hat C > 1.$ Let, for every $i\in[\ell]$, 
$$
C_{i}^- = ((1 - \delta)\hat C)^{(t-1)^{i-1}},\quad 
C_{i}^+ = ((1 +\delta)\hat C)^{t^{i-1}};
$$ 
$$
A_i = A_{i-1} \cap \{C_{i}^- x \le X_i \le C_{i}^+ x\}.
$$

Let us prove that 
\begin{equation}
\mathsf{P}(C_{i}^-  x \le X_i \le C_{i}^+ x | A_{i-1})  \geq 1- (1+o(1)) \frac{1}{\delta^2 C_{i-1}^- \hat C x},
\label{A_general_induction}
\end{equation}
uniformly over all $i\in\{2,3,\ldots,\ell\}$.
 
%For $i = 1$, we get $C_{1}^- = (1 - \delta) \hat C$, $C_{1}^+ = (1 + \delta) \hat C $, and (\ref{A_general_induction}) follows from (\ref{A_1st_step}). 

%Now suppose that $i\in[\ell-1]$ and (\ref{A_general_induction}) holds.  
Notice that edges between $V(G_2) \setminus \bigcup\limits_{j=1}^{i} U_j$ and $U_i$ still have independent Bernoulli distributions. Take an integer $a\in[ C_i^-  x, C_i^+ x]$. Since, on $A_i$, almost all vertices of $V(G_2)$ can be included in $U_{i+1}$ (i.e., $|U_1\sqcup\ldots\sqcup U_i|\leq \ln n(C_1^++\ldots+C_i^+)=o(n)$), we get
$$
\Exp(X_{i+1}|A_i \cap\{X_i=a\}) \sim \frac{n}{2} \sum\limits_{k = t-1}^{a} {a \choose k} p^k(1-p)^{a-k} \sim \frac{n}{2}\frac{a^{t-1}}{(t-1)!}p^{t-1};
$$
$$
\frac{\Var(X_{i+1}|A_i\cap\{X_i = a\})}{({\sf E}(X_{i+1}|A_i\cap\{X_i=a\}))^2} \sim \frac{1}{{\sf E}(X_{i+1}|A_i\cap\{X_i = a\})}.
$$

So, 
$$
\Exp(X_{i+1}|A_i) \le  (1 + o(1)) \frac{n}{2} \frac{(C_i^+ x)^{t-1}}{(t-1)!}p^{t-1} = (1 + o(1)) \frac{C_{i+1}^+}{(1 + \delta)^{t^{i-1}}\hat C^{t^{i-1}-1}} x.
$$ 
%$$
%= (1 + o(1)) (1+\delta)^{t^{(i-2)} (t-1)} \hat C^{t^{i-2} (t-1) + 1} x 
%= (1 + o(1)) \frac{C_{i+}}{1 + \delta} x.
%$$  

Similarly, 
$$
\Exp(X_{i+1}|A_i) \ge  (1 + o(1)) \frac{n}{2} \frac{(C_i^- x)^{t-1}}{(t-1)!}p^{t-1} = (1 + o(1))  C_{i+1}^- \hat C x.
$$ 
%$$ 
%= (1 + o(1)) (1-\delta)^{(t-1)^{(i-2)} (t-1)} \hat C^{(t-1)^{i-2} (t-1)+1} x = %$$ 
%$$ 
%=  (1 + o(1)) { C_{i-} \hat C}x \ge (1 + o(1)) { C_{i-}}x.
%$$

%Notice that
%Similarly  $\Exp(X_i|A_{i-1}, X_{i-1}=a) \ge  (1 + o(1)) (1 - \delta)^{t-1} \hat C^t x $.

Therefore, by Chebyshev's inequality, 
$$
 {\sf P}\biggl(\{X_{i+1}>C_{i+1}^+x\}\cup\{X_{i+1}<C_{i+1}^-x\}\biggl|A_i\biggr)\leq
$$
$$
 {\sf P}\biggl(|X_{i+1}-{\sf E}(X_{i+1}|A_i)>\delta{\sf E}(X_{i+1}|A_i)\biggl|A_i\biggr)\leq
\frac{\Var(X_{i+1}|A_i)}{(\delta{\sf E}(X_{i+1}|A_i))^2} \le \frac{1+o(1)}{\delta^2 C_{i+1}^- \hat C x}.
$$
This finishes the proof of (\ref{A_general_induction}). \\

Notice that (\ref{A_1st_step}) and (\ref{A_general_induction}) imply
$$
\mathsf{P}(\neg A_{\ell}) \le  o(1) + \sum\limits_{i = 1}^{\ell}\mathsf{P}(C_{i-}  x \le X_i \le C_{i+} x|A_{i-1}) \le
$$ 
$$ 
\le o(1) + \frac{1}{\delta^2\hat C x}+ \sum_{i=2}^{\ell} \frac{1}{\delta^2 C_i^- \hat C x} \le o(1) + \frac{1}{\delta^2 \hat C x} \sum\limits_{k =0}^{+\infty} \frac{1}{((1-\delta) \hat C)^k}\to 0, \text{ as } (1-\delta) \hat C > 1.
$$
But, on $A_{\ell}$,

%Then $P(C_{i-}  x \le X_i \le C_{i+} x|A_{i-1}) =  %P(C_{i-1-} \hat C (1-\delta) x \le X_i \le  C_{i-1+} x|A_{i-1}) = P((1-\delta)(1 - \delta)^{s-1} \hat C^s x \le X_2 \le (1+\delta) (1 + \delta)^{s-1} \hat C^s x|A_1) \le P(|X_2 - EX_2| \ge (\delta + o(1)) EX_2|A_2) 
%\le \frac{D(X_2|A_1)}{(\delta + o(1))^2 E(X_2|X_1=a))^2}\sim \frac{1}{\delta^2 (1 - \delta)^{s-1} \hat C^s x}$

%By Chebyshev's inequality, 
%$$
%\mathsf{P}(C_{i-1}^-  x \le X_{i+1} \le C_{i+1}^+ x|A_{i-1}) \le  (1+o(1)) \frac{1}{\delta^2 C_{i-} \hat C x}
%$$ 
%which proves our statement.

%On  $A_{\ell}$, 
$$
|V(F_0)| + \sum\limits_{i=1}^{\ell} |U_i| = x+\sum\limits_{i=1}^{\ell} X_i \ge x + \sum_{i=1}^{\ell} C_{i}^- x  = 
$$ 
$$
= x\left(1 + \sum _{i=1}^{\ell} ((1-\delta)\hat C)^{(t-1)^{i-1}}\right) \ge  ((1-\delta)\hat C)^{(t-1)^{\ell-1}}x\geq
$$
$$
\geq ((1-\delta)\hat C)^{\frac{(\ln \ln n) ^{2}}{t-1}}x = 
\exp\left[\frac{\ln ((1- \delta) \hat C)}{t-1}(\ln \ln n) ^{2} + \ln \ln n\right] \gg (\ln n)^{\frac{t-1}{t-2}}.
$$

%It is left to prove that $\mathsf{P}(A_{\ell})$ is small.

%$$\mathsf{P}(\neg A_{\ell}) \le  o(1) + \sum\limits_{i = 1}^{\ell}\mathsf{P}(C_{i-}  x \le X_i \le C_{i+} x|A_{i-1}) \le$$ $$ \le o(1) + \sum_{i \le \ell} \frac{1}{\delta^2 C_{i-} \hat C x} \le o(1) + \frac{1}{\delta^2 \hat C x} \sum\limits_{k =1}^{+\infty} \frac{1}{((1-\delta) \hat C)^k} \le  $$ $$ \le o(1) + \frac{1}{\delta^2 \hat C x} \to 0, \text{ as } (1-\delta) \hat C > 1.$$

So, whp there exists a saturating structure of size $y$.
% \frac{C^{s-1}}{(t-1)!} \frac{\ln n^{t-1}}{\omega_n^{(s-1)^2/(s-2)}} \frac{n^{-1}\ln n\omega_n^{t-1}}{\ln n^{s-1}} \sim \frac{C_1^{t-1}}{2(t-1)!}\frac{\ln n}{\omega_n^{\frac{s-1}{s-2}}} =  \frac{C_1^{t-2}}{2(t-1)!}x = \hat C x$

%$\ln x_0 = \frac{(2-s)(1 +o(1)) + (1-s)\sigma_n-1+\ln A}{(s-2)(1+o(1))}$

%$x_0 = const*(\frac{\ln n}{\gamma_n\omega_n})^{\frac{s-1}{s-2}}$

%$\ln p = \ln C + \frac{1}{t-1}\ln n - \frac{t-2}{t-1}\ln \ln n$

%$\prod_{i=t}^{x}{i-1 \choose t-1} = \frac{1}{(t-1)!^{x-t+1}} \frac{t!...(x-1)!}{1!...(x-t)!} = \frac{1}{(t-1)!^{x-t+1}} \frac{(x-t+1)!...(x-1)!}{1!...(t-1)!} \ge \frac{1}{A^{x-t+1}} (x-t+1)^{(x-t+1)(t-1)}$

\subsubsection{Saturating structure of size $n$}

In this section, we prove that, for every proper $S\subset[n]$ of size at least  $y$, in $G(n,p)$  there exists a vertex outide $S$ such that it is connected to at least $t-1$ vertices of $S$. Clearly, this observation finishes the proof of Claim~\ref{existence_saturating}.2.\\

%By the union bound, the probability of the contrary does not exceed  $${n \choose y} \left(1 - (1+o(1)){y \choose {t-1}} p^{t-1}\right)^{n-y} = $$ $$ =% O\left(\left(\frac{ne}{y}\right)^y \left(1 - (1+o(1)){(\frac{y^{t-1} }{(t-1)!})}p^{t-1}\right)^{n-y}\right)
% O\left(\exp \left \{ y(\ln n +1 - \ln y) - n(1+o(1)) \left(\frac{y^{t-1} }{(t-1)!}\right)p^{t-1} \right\} \right) = $$ $$ = O\left(\exp \{ y(\ln n + 1 - \ln y - (1+o(1)) \ln n\}\right)  = O(\exp\{-y\ln y\}) \to 0.$$

%Now suppose .

By the union bound, the probability that there exists a set $S\subset[n]$ of size $z\in[y,n/\ln n]$ such that every vertex outside $S$ has less than $t-1$ neighbors in $S$ is at most

%If $y \le z \le \frac{n}{\ln n}$, then the probability that there exists a set $S$ of size $z$ such that all vertices that do not belong to $S$ are connected to less than $t-1$ vertices of $S$ is at most 

$$
\sum_{z=y}^{\lfloor n/\ln n\rfloor}{n \choose z} \left(1 - {z \choose {t-1}} p^{t-1}\right)^{n-z}\leq 
$$
$$
\sum_{z=y}^{\lfloor n/\ln n\rfloor}
\exp \left[ z(\ln n +1 - \ln z) -  n (1+o(1)) \frac{z^{t-1} }{(t-1)!}p^{t-1} \right] = n\exp[-\Omega(z\ln n)])\to 0
$$
since $\frac{C^{t-1}}{(t-1)!}>1$.\\
 
Finally, the probability that there exists a proper $S\subset[n]$ of size $z>n/\ln n$ such that every vertex outside $S$ has less than $t-1$ neighbors in $S$ is at most

%If $z > \frac{n}{\ln n}$ this probability does not exceed 
 
%$${n \choose z} \left(1 - (1+o(1)){z \choose {t-1}} p^{t-1}\right)^{n-z} = $$

$$ % O\left(\left(\frac{ne}{y}\right)^y \left(1 - (1+o(1)){(\frac{y^{t-1} }{(t-1)!})}p^{t-1}\right)^{n-y}\right)
\sum_{z=\lceil n/\ln n\rceil}^{n-1} {n \choose n-z} \left(1 - {z \choose {t-1}} p^{t-1}\right)^{n-z}\leq
$$ 
$$
\sum_{z=\lceil n/\ln n\rceil}^{n-1} 
\exp \left[(n-z)\left(\ln n -  \frac{z^{t-1} }{(t-1)!}p^{t-1}\right) \right]=
n\exp\left[-\Omega\left(\frac{n^{t-2}}{(\ln n)^{2t-3}}\right)\right]\to 0.
$$

%Then the probability that there exists a set $S$ of size at least $y$ such that all vertices that do not belong to $S$ are connected to less than $t-1$ vertices of $S$ is at most 

%$$n e^{-y \ln y} + ne^{-n} \le e^{\ln n - y \ln y} + o(1) \to 0 .$$

%So, whp we can extend our saturating structure of size $y$ by one vertex at each step  (there exist $t-1$ edges to the current structure) and get a saturating structure of size $n$ and 

Claim~\ref{existence_saturating}.2 and Theorem~\ref{threshold_stars} follows.

%\section{Discussions}

\section{Acknowledgements}

This work is supported by RFBR, grant number 20-51-56017.

\end{document}